\def\version{17/11/2005 Version 11}

\documentclass[11pt,a4paper]{amsart}\usepackage{fullpage}
\usepackage{amssymb,amscd}
\usepackage[all]{xy}
\usepackage{mathrsfs}

\newtheorem{thm}{Theorem}[section]
\newtheorem{lem}[thm]{Lemma}
\newtheorem{prop}[thm]{Proposition}
\newtheorem{cor}[thm]{Corollary}

\theoremstyle{definition}
\newtheorem{rem}[thm]{Remark}
\newtheorem{rems}[thm]{Remarks}
\newtheorem{defn}[thm]{Definition}
\newtheorem{ex}[thm]{Example}

\newtheorem{cond}[thm]{Condition}
\newtheorem*{Cond*}{Condition}

\newtheorem*{Quest*}{Question}

\numberwithin{section}{part}
\numberwithin{equation}{section}
\numberwithin{figure}{section}

\def\ie{\emph{i.e.}}
\def\ds{\displaystyle}
\def\:{\colon}
\def\.{\cdot}
\def\o{\circ}
\def\<{\left\langle}
\def\>{\right\rangle}
\def\({\left(}
\def\){\right)}
\def\ph#1{\phantom{#1}}
\def\epsilon{\varepsilon}
\def\phi{\varphi}
\def\subset{\subseteq}

\def\leq{\leqslant}
\def\geq{\geqslant}
\def\lla{\longleftarrow}

\def\lra{\longrightarrow}
\def\Lra{\Longrightarrow}
\def\ra{\rightarrow}
\def\bar#1{\overline{#1}}
\def\hat#1{\widehat{#1}}
\def\tilde#1{\widetilde{#1}}
\def\iso{\cong}

\def\C{\mathbb{C}}

\def\F{\mathbb{F}}

\def\N{\mathbb{N}}
\def\Q{\mathbb{Q}}

\def\Z{\mathbb{Z}}
\def\Smash_#1{\wedge_{#1}}

\def\Times_#1{\ds\mathop{\times}_{#1}}
\def\oTimes_#1{\otimes_{#1}}
\def\oPlus_#1{\bigoplus_{#1}}
\def\ideal{\triangleleft}
\DeclareMathOperator{\inc}{inc}
\DeclareMathOperator{\mult}{mult}
\DeclareMathOperator{\tr}{tr}
\DeclareMathOperator{\Aut}{Aut}
\DeclareMathOperator{\Br}{Br}
\DeclareMathOperator{\cone}{cone}
\DeclareMathOperator{\End}{End}
\DeclareMathOperator{\Ext}{Ext}
\DeclareMathOperator{\Gal}{Gal}
\DeclareMathOperator{\Har}{Har}
\DeclareMathOperator{\Hom}{Hom}
\DeclareMathOperator{\Kumm}{Kum}
\DeclareMathOperator{\lcm}{lcm}
\DeclareMathOperator{\map}{map}
\DeclareMathOperator{\Map}{Map}
\DeclareMathOperator{\Pic}{Pic}

\DeclareMathOperator{\Tor}{Tor}
\def\id{\mathrm{id}}

\def\HG{\mathrm{H}\Gamma}

\DeclareMathOperator{\HH}{HH}

\DeclareMathOperator{\Der}{Der}

\def\AbGps{\mathbf{AbGps}}
\def\FinAbGps{\mathbf{FinAbGps}}
\DeclareMathOperator{\Real}{Real}

\title[Realizability of Galois extensions by ring spectra]
{Realizability of algebraic Galois extensions by strictly
commutative ring spectra}
\author{Andrew Baker \and Birgit Richter}
\address{Mathematics Department, University of Glasgow,
Glasgow G12 8QW, Scotland.}
\email{a.baker@maths.gla.ac.uk}
\urladdr{http://www.maths.gla.ac.uk/$\sim$ajb}
\address{Fachbereich Mathematik der Universit\"at Hamburg,
Bundesstrasse 55, 20146 Hamburg, Germany.}
\email{richter@math.uni-hamburg.de}
\urladdr{http://www.math.uni-hamburg.de/home/richter/}
\date{\version}
\thanks{
We would like to thank John Rognes, John Greenlees, Peter
Kropholler, Stefan Schwede and the referee for helpful comments.
We thank the Mathematics Departments of the Universities of Glasgow
and Oslo for providing us with stimulating environments to pursue
this work. \\
To appear in Transactions of the American Mathematical Society.
\texttt{math.AT/0406314}
}
\keywords{Commutative $S$-algebra, Galois extension, $\Gamma$-homology,
Kummer theory, Picard groups}
\subjclass[2000]{Primary 55P42, 55P43, 55S35; Secondary 55P91, 55P92, 13B05}

\begin{document}
\begin{abstract}
We discuss some of the basic ideas of Galois theory for commutative
$\mathbb{S}$-algebras originally formulated by John Rognes. We
restrict attention to the case of finite Galois groups and to global
Galois extensions.

We describe parts of the general framework developed by Rognes.
Central r\^{o}les are played by the notion of \emph{strong duality}
and a \emph{trace mapping} constructed by Greenlees and May in the
context of generalized Tate cohomology. We give some examples where
algebraic data on coefficient rings ensures strong topological
consequences. We consider the issue of passage from algebraic Galois
extensions to topological ones applying obstruction theories of
Robinson  and Goerss-Hopkins to produce topological models for
algebraic Galois extensions and the necessary morphisms of
commutative $\mathbb{S}$-algebras. Examples such as the complex
$K$-theory spectrum as a $KO$-algebra indicate that more exotic
phenomena occur in the topological setting. We show how in certain
cases topological abelian Galois extensions are classified by the
same Harrison groups as algebraic ones and this leads to computable
Harrison groups for such spectra. We end by proving an analogue of
Hilbert's theorem~90 for the units associated with a Galois
extension.
\end{abstract}

\maketitle

\section*{Introduction}

We discuss some ideas on Galois theory for commutative
$\mathbb{S}$-algebras, also known as brave new (commutative) rings,
originally formulated by John Rognes. We restrict ourselves to the
case of finite Galois groups, although there are versions for
profinite groups, group-like monoids, and stably dualizable
topological groups. The Galois extensions which we consider are
\emph{global} in Rognes' terminology, \ie, we work in the stable
homotopy category and not in dramatically localized versions of it.

We begin in Part~\ref{part:1} by describing the general framework,
first outlining the generalization of `classical' Galois theory of
field extensions to commutative rings and then some of the theory
developed by Rognes. In our account, central r\^{o}les are played by
the notion of \emph{strong duality} as discussed by Dold and
Puppe~\cite{D&P} and a certain \emph{trace mapping}  constructed by
Greenlees and May in the context of generalized Tate cohomology.
Both of these are topological manifestations of properties of
algebraic Galois extensions. We stress that this material is known
to Rognes and
is systematically described in~\cite{R:Opusmagnus}, and the presentation
reflects our
approach to understanding his results rather than being original.

In the heart of our paper, Part~\ref{part:2}, we consider the issue
(raised by Rognes) of passage from algebraic Galois extensions to
topological ones, \ie, given a commutative $\mathbb{S}$-algebra $A$
and a $G$-Galois extension $B_*$ of $A_* = \pi_*(A)$ we investigate
whether this extension can be realized by a Galois extension $B/A$
of $\mathbb{S}$-algebras. We focus on specific situations where we
can apply the recently developed obstruction theories of Robinson
and Goerss-Hopkins to produce topological models for algebraic
Galois extensions and the necessary morphisms of commutative
$\mathbb{S}$-algebras. In these situations, certain cohomological
obstructions vanish for purely algebraic reasons related to the
Galois theory of the coefficient rings of the spectra involved.
However, examples such as $KU$ as a $KO$-algebra (studied in
Part~\ref{part:1}) indicate that more exotic phenomena can occur in
the topological setting.

In Section~\ref{sec:TopAbelExtns&KummerThy} we investigate Kummer
extensions of commutative  $\mathbb{S}$-algebras and provide
criteria that allow to compare them with algebraic Kummer
extensions. In certain cases, for instance the $C_2$-extensions of
$KO[1/2]$,  topological abelian Galois extensions are classified by
the same Harrison groups as algebraic ones and this leads to
computable Harrison groups for such spectra. An important technical
input is provided by some results about invertible module spectra
and topological Picard groups proved in~\cite{AB&BR:Inv}. We end by
proving an analogue of Hilbert's theorem~90 for the units associated
with a Galois extension.

\part{Galois theory for commutative rings and ring spectra}
\label{part:1}

\section{Galois theory for commutative rings}\label{sec:Galois-CommRings}

We recall results on the Galois theory of commutative rings which
are mainly due to Chase, Harrison and Rosenberg~\cite{CHR} and also
described by Greither~\cite{Greither:LNM1534}. The notion of Galois
extensions of commutative rings was first developed by Auslander and
Goldman in~\cite{AG}; for further work on this
see~\cite{Ferrero&Paques}.
It is also possible to make sense of these ideas in the context where
$R$ is a commutative graded ring and $S$ is a commutative $R$-algebra,
and this will be important in the topological applications.

Let $R$ be a commutative ring and $S$ a commutative $R$-algebra. Suppose
that $G\leq\Aut(S/R)$, the group of all $R$-algebra automorphisms of~$S$;
we will indicate the (left) action of $\gamma\in G$ on $s\in S$ by writing
${}^\gamma s$. We give the product $S$-algebra
\[
\prod_{\gamma\in G}S=\{(s_\gamma)_{\gamma\in G}:s_\gamma\in S\}
\]
the left $G$-action for which
\[
\alpha\.(s_\gamma)_{\gamma\in G}=(s_{\gamma\alpha})_{\gamma\in G}
\quad(\alpha\in G).
\]
We also have the $S$-algebra of functions $f\:G\lra S$ which has the
left $G$-action
\[
(\alpha\.f)(\gamma)=f(\alpha^{-1}\gamma)\quad(\alpha\in G).
\]
There is a $G$-equivariant isomorphism of $S$-algebras
\begin{equation}\label{eqn:IsomMap-Prod}
\Map(G,S)\iso\prod_{\gamma\in G}S;
\quad
f\longleftrightarrow\(f(\gamma^{-1})\)_{\gamma\in G}
\end{equation}
and a unique $S$-algebra homomorphism
\begin{equation}\label{eqn:Theta}
\Theta\:S\oTimes_RS\lra\prod_{\gamma\in G}S;
\quad
\Theta(u\otimes v)=(u{}^{\gamma}v)_{\gamma\in G}\quad(u,v\in S),
\end{equation}
which induces a $G$-equivariant homomorphism of $S$-algebras
\[
\Theta'\:S\oTimes_RS\lra\Map(G,S)
\]
that will be used without further remark.

There is an isomorphism of $S[G]$-modules
\begin{equation}\label{eqn:groupring}
\Upsilon: S[G] \lra \Map(G,S),
\end{equation}
sending a generator $\gamma$ to the Kronecker function
$\delta_\gamma$. Here, we use the left $S$-linear $G$-action on the
group ring $S[G]$.

We denote the \emph{twisted group ring of $S$ with $G$} by $S\sharp
G$; this is the free $S$-module on $G$ with product given by
\[
(s\alpha)(t\beta)=(s\,{}^\alpha t)(\alpha\beta)
\quad
(s,t\in S,\;\alpha,\beta\in G).
\]
Then $S\sharp G$ is an $R$-algebra and there is an $R$-algebra homomorphism
$j\:S\sharp G\lra\End_R(S)$ induced from the actions of $S$ and $G$ on $S$
by $R$-module homomorphisms.
\begin{defn}\label{defn:GaloisExtn}
$S/R$ is a \emph{$G$-Galois extension} if it satisfies the following
conditions: \\[3pt]
(G-1) $S^G=R$; \\[3pt]
(G-2) $\Theta\:S\oTimes_RS\lra\prod_{\gamma\in G}S$ is an isomorphism.
\end{defn}
\begin{rem}\label{rem:GaloisExtn-1}
(a) Condition (G-2) ensures that $S$ is unramified with respect to $R$.
For instance, if $C_2$ is a cyclic group of order~$2$, then $\Z\lra\Z[i]$
is \emph{not} a $C_2$-Galois extension: because of ramification at the
prime~$2$, the map $\Theta$ is not surjective. After inverting~$2$ we
find that $\Z[1/2]\lra\Z[1/2,i]$ is a $C_2$-Galois extension.
See Example~\ref{ex:Cyclotomic} for related phenomena. \\
(b) (G-2) can also be replaced by the requirement that
$\Theta'\:S\oTimes_RS\lra\Map(G,S)$ is an isomorphism.
\end{rem}
\begin{rem}\label{rem:GaloisExtn-2}
Later we will consider $G$-Galois extensions of graded commutative
rings. By these we mean extensions of graded rings $R_*\lra S_*$
together with an action of $G\leq\Aut_{R_*}(S_*)$ such that the
conditions of Definition~\ref{defn:GaloisExtn} are satisfied. Note
that in these cases the $G$-action preserves the grading.
\end{rem}
\begin{thm}\label{thm:GaloisExtn-Eqces}
Let $R$ be a commutative ring and $S$ be a commutative  $R$-algebra
with $G \leq \Aut_R(S)$ and assume that $S^G = R$, then the
following conditions are equivalent.
\begin{enumerate}
\item[(a)]
$S/R$ is a $G$-Galois extension.
\item[(b)]
$\Theta\:S\oTimes_RS\lra\prod_{\gamma\in G}S$ is an epimorphism.
\item[(c)]
There are finite sequences $u_1,\ldots,u_n,v_1,\ldots,v_n\in S$
for which
\[
\sum_{i=1}^nu_i{}^\gamma v_i=
\begin{cases}
1& \text{\rm if $\gamma=1$}, \\
0& \text{\rm otherwise}.
\end{cases}
\]
\item[(d)]
$S$ is a finitely generated projective $R$-module and
$j\:S\sharp G\lra\End_R(S)$ is an isomorphism.
\end{enumerate}
\end{thm}

We define the \emph{trace of $S/R$} to be the $R$-module homomorphism
\[
\tr_{S/R}\:S\lra R;\quad \tr_{S/R}(s)=\sum_{\gamma\in G}{}^\gamma s.
\]
We will write $\tr=\tr_{S/R}$ when no ambiguity is likely to result.

Recall that a ring is \emph{connected} if its only idempotents are $0$
and $1$. Also recall that given a left $R[G]$-module $M$, $\Hom_R(M,R)$
can be viewed as a left $R[G]$-module with the contragredient $G$-action
determined by
\[
(\gamma\.f)(m)=f(\gamma^{-1}m)\quad(f\in\Hom_R(M,R),\;\gamma\in G,\;m\in M).
\]
\begin{thm}\label{thm:GaloisExtn-Properties}
For a $G$-Galois extension $S/R$, the following hold.
\begin{enumerate}
\item[(a)]
$S$ is faithfully flat over $R$.
\item[(b)]
$\tr\:S\lra R$ is an epimorphism.
\item[(c)]
The unit $R\lra S$ is a split monomorphism of $R$-modules.
\item[(d)]
If $R$ and $S$ are both connected, then $\Aut(S/R)=G$.
\item[(e)]
For any commutative $R$-algebra $T$, $T\oTimes_RS/T$ is a $G$-Galois
extension.
\item[(f)]
$S$ is a finitely generated projective invertible $R[G]$-module, hence
it is of constant rank~$1$. Furthermore, $S$ is self-dual, \ie,
$S^*=\Hom_R(S,R)\iso S$ as $R[G]$-modules.
\end{enumerate}
\end{thm}
\begin{rem}\label{rem:GaloisExtn-GneqAut}
(i) The proof of (f) makes use of the trace pairing
\[
S\oTimes_R S\xrightarrow{\mult}S\xrightarrow{\tr_{S/R}}R
\]
to establish the self-duality and projectivity (see~\cite[theorem~4.2]{CHR}).
In Section~\ref{sec:Duality} we will discuss a more abstract setting in
which such duality occurs; in particular, the self-duality of a Galois
extension forces it to be self-dual in the categorical sense we will
describe (this is a special case of~\cite[lemma~2.9]{JPM:PicardGps} where
the relationship between these notations is also studied in detail). \\
(ii) When $R$ contains non-trivial idempotents,
\ref{thm:GaloisExtn-Properties}(d) need not be true. For example, if
$G$ has order~$2$ and $S=R\times R$ is the trivial $G$-extension,
then given any non-trivial idempotent $e\in R$, the map
\[
\phi\:R\times R\lra R\times R;
\quad
\phi(x,y)=(xe+y(1-e),x(1-e)+ye)
\]
is an $R$-algebra isomorphism which is not induced by an element
of~$G$.
\end{rem}

Condition~(f) has an important group cohomological consequence.
First recall the following well-known observation about group
cohomology (for example, see~\cite[example~6.1.2]{Weibel}).
\begin{lem}\label{lem:Cohom-Z->A}
For any ring~$R$ and $R[G]$-module $M$,
\[
\mathrm{H}^*(G;M)=\Ext^*_{\Z G}(\Z,M)\iso\Ext^*_{R[G]}(R,M).
\]
\end{lem}

Then we have
\begin{prop}\label{prop:G-Cohom-S/R}
Let $S/R$ be a $G$-Galois extension. Then
\[
\mathrm{H}^*(G;S)=\mathrm{H}^0(G;S)=R.
\]
\end{prop}
\begin{proof}
By Lemma~\ref{lem:Cohom-Z->A},
\[
\mathrm{H}^*(G;S)=\Ext_{R[G]}^*(R,S).
\]
Recall from Theorem~\ref{thm:GaloisExtn-Properties}(f) that $S$ is
finitely generated, self-dual and projective as an $R[G]$-module and
it is also finitely generated and projective as an $R$-module. As $S$
is a retract of a finitely generated free $R[G]$-module, it suffices
to prove the claim for $R[G]$. As $R[G] \cong \Map(G,R)$, an
adjunction argument proves the claim.
\end{proof}

For further results, as for instance the fundamental theorem of
Galois theory in this context see
\cite{CHR,Greither:LNM1534,Ferrero&Paques}. We
follow~\cite{Greither:LNM1534} in making Definitions~\ref{defn:GAL}
and~\ref{defn:Harrison} below.
\begin{defn}\label{defn:GAL}
For a commutative ring $R$ and a finite group $G$, let $\Gal(R,G)$
denote the category of $G$-Galois extensions of $R$ with morphisms
the $R$-algebra homomorphisms commuting with the actions of~$G$.
\end{defn}
\begin{prop}\label{prop:GAL-Groupoid}
If $S/R$ and $T/R$ are two $G$-Galois extensions and $\phi\:S\lra T$
is an $R$-algebra homomorphism commuting with the actions of~$G$,
then $\phi$ is an isomorphism. Hence $\Gal(R,G)$ is a {\rm(}large{\rm)}
groupoid.
\end{prop}
\begin{proof}
The proof of~\cite[proposition~0.1.12]{Greither:LNM1534} only applies
when~$R$ has no non-trivial idempotents, so for completeness we prove
the general case.

First note that by parts (a) and (e) of
Theorem~\ref{thm:GaloisExtn-Properties}, it suffices to replace
$S/R$ and $T/R$ by $(S\oTimes_RT\oTimes_RS)/(S\oTimes_RT)$ and
$(S\oTimes_RT\oTimes_RT)/(S\oTimes_RT)$, then note that
\[
S\oTimes_RT\oTimes_RS\iso\prod_{\gamma\in G}S\oTimes_RT
                                         \iso S\oTimes_RT\oTimes_RT.
\]
Thus we might as well assume that $S=T=\prod_{G}R$ is the trivial
$G$-Galois extension. For each $\alpha\in G$, there is an idempotent
$e_\alpha=(\delta_{\alpha\,\gamma})\in S$.

Now let
\[
e_1' = \phi(e_1) =\sum_{\gamma\in G}t_\gamma e_\gamma
\]
with $t_\gamma \in R$. Then as $\phi$ commutes with the action of
elements of $G$, for each $\alpha\in G$ we have
\[
e_\alpha' = \phi(e_\alpha) = \phi(\alpha\.e_1) = \alpha \phi(e_1)
=\sum_{\gamma\in G}t_{\alpha^{-1}\gamma}e_\gamma.
\]
Thus the $e'_\alpha$ form a complete set of orthogonal idempotents
in $S$. Notice that the $e_\alpha$ also form a basis for the free
$R$-module $S$. The equation $e'_1e'_1 = e'_1$ then shows that
$t_\gamma^2 = t_\gamma$ and
\[
\sum_\beta t_\beta = \sum_\beta t_\beta \sum_\gamma e_\gamma =
\sum_{\alpha, \gamma} t_{\alpha^{-1}\gamma}e_\gamma = \sum_\alpha
e'_\alpha = 1
\]
proves that the $t_\alpha$ also form a complete set of orthogonal
idempotents in~$R$.

Now for elements $s\in S$ and $x_\alpha\in R$, consider the equation
\begin{equation}\label{eqn:GAL-GroupoidProof}
\sum_{\alpha\in G}x_\alpha e'_\alpha=s,
\end{equation}
which is equivalent to
\[
\sum_{\alpha\in G}\sum_{\gamma\in G}x_\alpha t_{\alpha^{-1}\gamma}e_\gamma=s.
\]
Multiplying by $e_\beta$ we obtain $x_\alpha
t_{\alpha^{-1}\beta}e_\beta=r_\beta e_\beta$, where $r_\beta\in R$
is the unique element for which $se_\beta=r_\beta e_\beta$. Thus for
$\alpha,\beta\in G$ we have $x_\alpha t_{\alpha^{-1}\beta}=r_\beta$.
Multiplying by $t_{\alpha^{-1}\beta}$ now gives $x_\alpha
t_{\alpha^{-1}\beta}=r_\beta t_{\alpha^{-1}\beta}$. Summing over
$\beta\in G$ we now obtain $x_\alpha=\sum_{\beta\in G}r_\beta
t_{\alpha^{-1}\beta}$,  since $\sum_{\beta\in
G}t_{\alpha^{-1}\beta}=\sum_{\gamma\in G}t_\gamma=1$.
Thus~\eqref{eqn:GAL-GroupoidProof} has the unique solution given by
this formula. Hence $\phi$ is an isomorphism.
\end{proof}

Because of the last result, we may define an equivalence relation
$\,\sim\,$ on the objects of $\Gal(R,G)$ by requiring that $S/R\sim T/R$
if and only if there is a morphism $\phi\:S\lra T$ in $\Gal(R,G)$.
The equivalence classes are then the isomorphism classes of $G$-Galois
extensions of~$R$.
\begin{defn}[{See~\cite{Harrison:AbExtns} and~\cite[3.2]{Greither:LNM1534}}]
\label{defn:Harrison}
The \emph{Harrison set} $\Har(R,G)$ is the set of isomorphism classes
of $G$-Galois extensions of~$R$. When~$G$ is abelian, this is naturally
an abelian group often called the \emph{Harrison group}.
\end{defn}

There are some useful properties of this construction that will be
required later. Details can be found in~\cite[theorem~4]{Harrison:AbExtns}
or~\cite[theorems 3.2,3.3 and 3.5]{Greither:LNM1534}
or supplied by the reader.
\begin{prop}\label{prop:harrison}
Let $R$ be a commutative ring.
\begin{enumerate}
\item[(a)]
$\Har(R,{-})$ defines a covariant functor
\[
\Har(R,{-})\:\FinAbGps\rightsquigarrow\AbGps
\]
from finite abelian groups to abelian groups.
\item[(b)]
$\Har(R,{-})$ is left exact and is pro-representable.
\item[(c)]
$\Har(R,{-})$ preserves products, \ie, for finite abelian groups~$G$
and~$H$ there is a natural isomorphism
\[
\Har(R,G\times H)\iso\Har(R,G)\times\Har(R,H).
\]
\end{enumerate}
\end{prop}

Of course, part (c) implies that $\Har(R,{-})$ is determined by its
values on cyclic groups.

%

\section{Abelian extensions and Kummer theory}\label{sec:Abel&Kummer}

In \cite{Harrison:AbExtns,Greither:Kummer,Greither:LNM1534} a theory
of abelian extensions of commutative rings was described, including
an analogue of Kummer theory. We will describe this algebraic theory
and in Section~\ref{sec:TopAbelExtns&KummerThy} a topological analogue
will be introduced. Our goal is to show how $\Har(R,G)$ can be determined
under certain conditions, making use of Proposition~\ref{prop:harrison},
which reduces the problem to the case of cyclic groups.

Let $R$ be a commutative ring containing $1/n$ and a primitive $n$-th
root of unity $\zeta$. For a unit $u\in R^\times$, we set
\[
R(n;u)=R[x]/(x^n-u),
\]
where $x$ is an indeterminate. We will write $\bar z$ for the coset
$z+(x^n-u)\in R[x]/(x^n-u)$. Of course, for any $t\in R^\times$
there is a canonical $R$-algebra isomorphism
\[
R(n;t^nu)\xrightarrow{\iso}R(n;u); \quad \bar{tx} \mapsto \bar{x}.
\]
Notice that $R(n;u)/R$ is a $C_n$-Galois extension, where the action
of the generator $\gamma_n\in C_n$ is given by
\[
\gamma_n\.\bar{x}=\bar{\zeta x}.
\]

The set $\Kumm_n(R)$ of $R$-algebra isomorphism classes of such $R(n;u)$
is an abelian group with product on isomorphism classes given by
\[
[R(n;u)][R(n;v)]=[R(n;uv)]
\]
and whose unit is the class $[R(n;1)]$, where
\[
R(n;1)=R[x]/(x^n-1)=\prod_{i=1}^nR[x]/(x-\zeta^i)=\prod_{\gamma\in G}R
\]
is the trivial $G$-Galois extension. Of course, $\Kumm_n(R)\leq\Har(R,C_n)$.
In fact there is an isomorphism of groups
\begin{equation}\label{eqn:Kummer}
R^\times/(R^\times)^n\xrightarrow{\iso}\Kumm_n(R);
\quad
u(R^\times)^n\longmapsto[R(n;u)].
\end{equation}

Now let $S/R$ be a $C_n$-Galois extension. For $k=0,1,\ldots,n-1$, let
\[
S^{(k)}=\{s\in S:{}^{\gamma_n}s=\zeta^ks\}\subset S.
\]
Then each $S^{(k)}$ is an $R$-submodule of $S$ and is a summand,
hence it is finitely generated projective. Furthermore, the product
in~$S$ gives rise to isomorphisms $ S^{(k)}\oTimes_R S^{(\ell)}\lra
S^{(k+\ell)}$. In particular we obtain
\[
\overset{n}{\overbrace{S^{(1)}\oTimes_R\cdots\oTimes_R S^{(1)}}}
\xrightarrow{\iso}S^{(1)}\oTimes_RS^{(n-1)}
\xrightarrow{\iso}S^{(0)}=R.
\]
This shows that $S^{(1)}$ is an invertible $R$-module which represents
an element $[S^{(1)}]$ of the Picard group $\Pic(R)$ whose order is
a divisor of~$n$; we write
\[
\Pic(R)[n]=\{P\in\Pic(R):P^n=1\}.
\]
Thus there is a group homomorphism
\[
\Har(R,C_n)\lra\Pic(R)[n];\quad[S]\longmapsto[S^{(1)}].
\]

Now from~\cite{Greither:Kummer,Greither:LNM1534} we have
\begin{prop}\label{prop:KummSeq}
There is an exact sequence of abelian groups
\[
1\ra R^\times/(R^\times)^n\lra\Har(R,C_n)\lra\Pic(R)[n]\ra1.
\]
\end{prop}

There is a generalization of these ideas to the case where~$G$
is any finite abelian group and $R$ contains $1/|G|$ as well
as a primitive $d$-th root of unity $\zeta$ where
$d=\lcm\{|\gamma|:\gamma\in G\}$ is the exponent of $G$.
Presumably the following is known to experts (it is hinted at
in~\cite{Greither:Kummer,Greither:LNM1534}) but we give details
since we know of no convenient reference.

First note that as an $R[G]$-module, the group ring $R[G]$ has
a decomposition
\begin{equation}\label{eqn:R[G]-Decomp}
R[G]=\bigoplus_{\chi}R[G]e_\chi,
\end{equation}
where the sum is over the characters
$\chi\in\Hom(G,\<\zeta\>)=\Hom(G,R^\times)$.
This decomposition is effected by the orthogonal idempotents
\begin{equation}\label{eqn:R[G]-Idempot}
e_\chi=\frac{1}{|G|}\sum_{\gamma\in G}\chi(\gamma^{-1})\gamma\in R[G]
\end{equation}
which decompose~$1$. It is easily seen by direct calculation that
the $R$-module $R[G]e_\chi=e_\chi R[G]$ is free of rank~$1$.

Now for a $G$-Galois extension $S/R$ as above there is a decomposition
of $R[G]$-modules
\[
S=\bigoplus_{\chi}S(\chi),
\]
where $S(\chi)=e_\chi S$.

\begin{lem}\label{lem:KummGenMult}
For characters\/ $\chi_1,\chi_2$, the multiplication map
$S(\chi_1)\oTimes_{R}S(\chi_2)\ra S(\chi_1\chi_2)$ is an
isomorphism. Hence for each character $\chi$, $S(\chi)$ is an
invertible $R$-module.
\end{lem}
\begin{proof}
This is similar to the proof for the case of a cyclic group. The
invertibility comes about because each character $\chi$ has an
inverse character $\bar{\chi}$ defined by
$\bar{\chi}(\gamma)=\chi(\gamma)^{-1}$.
\end{proof}

The \emph{character group of\/ $G$} is the abelian group
$G^\o=\Hom(G,\Q/\Z)$. Then $G^\o$ is finite of order $|G|$ and
$G^\o\iso\Hom(G,\<\zeta\>)$. Now set
\[
\Pic(R,G)=\Hom(G^\o,\Pic(R))=\Hom(G^\o,\Pic(R)[d]).
\]
In order to give an estimate of $\Har(R,G)$ we state the following
result. As we will not use this result later on, we refrain from
giving a proof.
\begin{prop}\label{prop:KummSeq-Gen}
When $G$ is abelian, there is a natural exact sequence of abelian
groups
\[
0\ra\mathrm{H}^2(G^\o,R^\times)\lra\Har(R,G)\lra\Pic(R,G)\ra0.
\]
\end{prop}
We now briefly discuss the graded version of Kummer theory. For a
graded commutative ring $R_*$, we have to distinguish between the
cases where the characteristic is two and the general case. In the
former case the grading is easily dealt with, so we concentrate on
cases where two is not zero. Then the units of $R_*$ are in even
degrees only.

If we want to build the analogue of $R(n;u)$ for a graded ring $R_*$,
we need to assume that the degree of $u \in R_*^{\times}$ is divisible
by $2n$. Let $R_{r*}$ be the subring of $R_*$ of elements in degrees
divisible by $r$. We can still identify $\Kumm_n(R_*)$ with a quotient
of units
$$
R_{2n*}^\times/(R_{2*}^\times)^n \cong \Kumm_n(R_*).
$$

We still obtain an eigenspace decomposition of every element in
$\Har(R_*,C_n)$ and therefore every $C_n$-Galois extension gives
rise to a graded invertible module over $R_*$ of order dividing $n$.
Therefore we get a left-exact sequence
$$
1 \ra  R_{2n*}^\times/(R_{2*}^\times)^n \lra \Har(R_*,C_n)
                                                   \lra \Pic(R_*)[n],
$$
but not every element in $\Pic(R_*)[n]$ has to be in the image of
$\Har(R_*,C_n)$. For instance if $R_*$ is periodic such that
$R_{*+n} = R_*$ with $n$ even, then $\Sigma R_*$ is in $\Pic(R_*)[n]$
but cannot come from a $C_n$-extension. However, if we restrict attention
to elements in the Picard group which are concentrated in even degrees
and of order dividing $n$, the construction given
in~\cite[p.22]{Greither:LNM1534} ensures that such elements stem from
the Harrison group.

\section{Duality in a symmetric monoidal category}\label{sec:Duality}

The Galois theory for commutative rings of Section~\ref{sec:Galois-CommRings}
has some crucial aspects which can be generalized to the context of
symmetric monoidal categories. At the heart of this are the notions
of \emph{strong duality} and \emph{self-duality} both of which are
visible in the above account. The appropriate notions are described
in detail in~\cite{D&P} and some aspects appear in~\cite{HPS,LMS}.
We give an account based on~\cite{D&P} but with some modifications
of notation.

Let $\mathscr{C}$ be a closed symmetric monoidal category with
multiplication $\boxtimes$, twist map $\tau$ and unit $I$. We denote
the internal homomorphism object on $X,Y \in \mathscr{C}$ by $F(X,Y)$
and write $DX$ for $F(X,I)$.
%
For every  $X \in \mathscr{C}$ we have a canonical \emph{evaluation
morphism} $\epsilon=\epsilon_X\: DX \boxtimes X \lra I$ which corresponds
to $\id_X\in\mathscr{C}(X,X)\iso\mathscr{C}(DX \boxtimes X,I)$. An object
$X$ is \emph{weakly self-dual} if $X$ is isomorphic to $DX$.
%
%
There is always a map $\delta=\delta_X\:X\lra DDX$ which corresponds to
\[
X \boxtimes DX \xrightarrow[\iso]{\tau} DX \boxtimes X
\xrightarrow{\epsilon} I.
\]
\begin{defn}\label{defn:Reflexive}
An object $X$ is \emph{reflexive} if $\delta_X$ is an isomorphism.
\end{defn}

We may define $\mu=\mu_{XY}\:DX \boxtimes DY \lra D (Y \boxtimes X)$
corresponding to the composite
\[
DX \boxtimes DY \boxtimes Y \boxtimes X
\xrightarrow{\id\boxtimes\epsilon_Y\boxtimes\id}
DX \boxtimes I \boxtimes X \iso DX  \boxtimes X
\xrightarrow{\epsilon_X}I.
\]
\begin{defn}\label{defn:StrongDual}
An object $X$ is \emph{strongly dualizable} if it is reflexive and
$\mu_{X\,DX}$ is an isomorphism. This condition is equivalent to the
requirement that the composition
\[
DX \boxtimes X \xrightarrow{\id\boxtimes\delta_X}
DX \boxtimes DDX \xrightarrow{\mu}
D(DX \boxtimes X)
\]
be an isomorphism and this means that $DX \boxtimes X$ is
canonically weakly self-dual.

If $X$ is strongly dualizable and weakly self-dual, we call $X$
\emph{strongly self-dual}.

\end{defn}

When $X$ is strongly dualizable, the \emph{coevaluation}
$\eta=\eta_X\:I\lra X\boxtimes DX$ is the composite
\[
I = DI \xrightarrow{D\epsilon} D(DX \boxtimes X)
\xrightarrow[\iso]{\mu^{-1}}DX \boxtimes DDX
\xrightarrow[\iso]{\id\boxtimes\delta^{-1}}DX \boxtimes X
\xrightarrow[\iso]{\tau} X \boxtimes DX.
\]

The following result taken from~\cite[theorem~1.3]{D&P} summarizes
the main properties of duality.
\begin{thm}\label{thm:Duality}
Let $X$ be an object of\/ $\mathscr{C}$ and $\epsilon\:DX \boxtimes X\lra I$
the evaluation. Then the following conditions are equivalent.
\begin{enumerate}
\item[(a)]
$X$ is  strongly dualizable.
\item[(b)]
There is a morphism $\eta\:I\lra X \boxtimes DX$ for which the
following compositions are the identity morphisms $\id_X$ and
$\id_{DX}$ respectively:
\begin{subequations}\label{eqn:Duality}
\begin{align}
\id_X \: X \iso I\boxtimes X \xrightarrow{\eta\boxtimes\id_X}
& X \boxtimes DX \boxtimes X
\xrightarrow{\id_X \boxtimes\epsilon} X \boxtimes I\iso X,
\label{eqn:Duality-P}\\
\id_{DX} \: DX \iso DX \boxtimes I \xrightarrow{\id_{DX} \boxtimes\eta}
& DX \boxtimes X \boxtimes DX
\xrightarrow{\epsilon\boxtimes\id_{DX}}I \boxtimes DX \iso DX.
\label{eqn:Duality-Q}
\end{align}
\end{subequations}
\item[(c)]
For every pair of objects $U$ and $V$, the map $
\phi_{UV}\:\mathscr{C}(U,V \boxtimes DX) \ra \mathscr{C}(U \boxtimes
X, V)$, which sends $f\:U \lra V \boxtimes DX$ to the composite
\[
U \boxtimes X \xrightarrow{f\boxtimes\id_X} V \boxtimes DX \boxtimes
X \xrightarrow{\id_V \boxtimes\epsilon} V \boxtimes I\iso V,
\]
is a bijection.
\end{enumerate}
Furthermore, if these conditions are satisfied then the morphism
$\eta$ of {\rm(b)} is necessarily the coevaluation and the bijection
$\phi_{UV}$ of {\rm(c)} sends $\eta$ to the composition $I\boxtimes
X \iso X \xrightarrow{\id_X} X$.
\end{thm}


We end this section with a useful result mentioned in the proof
of~\cite[theorem~XVI~7.4]{JPM:EqtHtpyCohThy}.
\begin{prop}\label{prop:StrongDual-Retract}
Let $\mathscr{C}$ be closed and let $X$ be strongly dualizable
in $\mathscr{C}$. Suppose that $R$ is a retract of $X$ with
maps $j\:R\lra X$ and $r\:X\lra R$ satisfying $rj=\id_R\:R\lra R$.
Then $R$ is strongly dualizable with the natural evaluation
$\epsilon'\: DR \boxtimes R \lra I$ and coevaluation given by the
composite
\[
\eta'\:I\xrightarrow{\eta} X \boxtimes DX
\xrightarrow{r \boxtimes j^*}R\boxtimes F(R,I),
\]
where $j^* \: DX \lra DR$ is dual to $j\:R\lra X$.
\end{prop}

\section{Brave new Galois extensions}\label{sec:GaloisExtns}

The notion of a Galois extension in the context of the commutative
$\mathbb{S}$-algebras of~\cite{EKMM} was introduced by John Rognes. We
restrict attention to the case of finite Galois groups and all
Galois extensions which we consider are `global', \ie, we work in
the unlocalized setting.

Given a commutative $\mathbb{S}$-algebra $A$, we will work in the
categories of $A$-modules $\mathscr{M}_{A}$ and its derived category
$\mathscr{D}_{A}$. These categories are complete and symmetric
monoidal under the smash product $\Smash_{A}$. In $\mathscr{D}_{A}$,
an $A$-module $L$ has as a weak dual $D_AL=F_A(L,A)$.

Then $L$ is strongly dualizable if for every $A$-module $M$,
$F_A(L,M)\sim F_A(L,A)\Smash_A M$, while $L$ is strongly self-dual
if in addition $F_A(L,A)\sim L$.

The following useful result on strongly dualizable objects in $\mathscr{D}_{A}$
is taken from \cite[\S 2]{HPS} (see also
\cite{JPM:EqtHtpyCohThy} and~\cite{JPM:PicardGps}).
\begin{prop}\label{prop:StrongDual-Cell}
Let $X$ be an $A$-module. Then $X$ is strongly dualizable in $\mathscr{D}_{A}$
if and only if it is weakly equivalent to a retract of a finite cell
$A$-module.
\end{prop}
%
%

Our next definition is of course suggested by the algebraic notion of
faithful flatness.
\begin{defn}\label{defn:Faithful}
An $A$-module $N$ is \emph{faithful} (as an $A$-module) if whenever~$M$
is an $A$-module for which $N\Smash_AM\sim*$, then $M\sim*$.
\end{defn}
\begin{rem}\label{rem:Faithful}
If $N$ is faithful, then the homology theory $N^A_*(-)=\pi_*(N\wedge_A-)$
detects weak equivalences since a morphism of $A$-modules $f\:M\lra M'$
is a weak equivalence if and only if the induced homomorphism
$f_*\:N^A_*M\lra N^A_*M'$ is an isomorphism.
\end{rem}

We can now give the key definition of a Galois extension essentially
due to Rognes~\cite{R:Opusmagnus}.
\begin{defn}\label{defn:BraveNewGaloisExtn-Rognes}
Let $A$ be a commutative $\mathbb{S}$-algebra and let $B$ be a commutative
cofibrant $A$-algebra. Let $G$ be a finite (discrete) group and suppose
that there is an action of $G$ on $B$ by commutative $A$-algebra morphisms.
Then $B/A$ is a \emph{weak $G$-Galois extension} if it satisfies the
following two conditions. \\[3pt]
(BNG-1) The natural map $A\lra B^{\mathrm{h}G}=F(E  G_+,B)^G$
is a weak equivalence of $A$-algebras; \\[3pt]
(BNG-2) There is a natural equivalence of $B$-algebras
$\Theta\:B\Smash_A B\xrightarrow{\sim}F(G_+,B)$ induced from the action
of $G$ on the right hand factor of~$B$. \\[3pt]
$B/A$ is a \emph{$G$-Galois extension} if it also satisfies \\[3pt]
(BNG-3) $B$ is faithful as an $A$-module.
\end{defn}

In fact Rognes does not insist on (BNG-3) but calls extensions satisfying
(BNG-1) and (BNG-2) Galois extensions and adds faithfulness as a requirement
whenever needed. So far, there are no known examples of Galois extensions
which are not faithful, though.

In (BNG-2), we use the topological analogue of the map $\Theta$
from~\eqref{eqn:Theta}. We also consider the maps $\tilde\gamma$
defined in~\eqref{eqn:EtaleMap-factor}. These have product
\begin{equation}\label{eqn:EtaleMap}
B\Smash_AB\lra\prod_{\gamma\in G}B.
\end{equation}
The following base-change results can be found in \cite[\S
7]{R:Opusmagnus}.

\begin{prop}\label{prop:FFbaseChange}
Let $A$ be a commutative $\mathbb{S}$-algebra and let $A\lra B$ and
$A\lra C$ be maps of commutative $\mathbb{S}$-algebras.
\begin{enumerate}
\item[(a)]
$C\Smash_{A}B$ admits a canonical commutative $C$-algebra structure.
\item[(b)]
If\/ $G$ acts on $B$ by $A$-algebra morphisms, then there is a
canonical extension of the action of\/ $G$ on $B$ to one by
$C$-algebra morphisms on $C\Smash_{A}B$.
\item[(c)]
If\/ $G$ acts on $B$ by $A$-algebra morphisms and $C$ is strongly
self-dual in $\mathscr{D}_{A}$ and $C$ is faithful as an $A$-module,
then $B/A$ is a $G$-Galois extension if and only if $C\Smash_{A}B/C$
is a $G$-Galois extension.
\end{enumerate}
\end{prop}

Here are some examples. Proofs that these are actually Galois
extensions can be found in~\cite{R:Opusmagnus}.

\begin{ex}\label{ex:TrivialGaloisExtn}
For a commutative $\mathbb{S}$-algebra $A$ and finite group $G$,
the morphism
\[
A\lra F(G_+,A)\iso\prod_{\gamma\in G}A
\]
induced from the trivial action of $G$ on $A$ is the \emph{trivial
$G$-Galois extension}.
\end{ex}
\begin{ex}\label{ex:HS/HR}
Let $R\lra S$ be a $G$-Galois extension in the algebraic sense of
Section~\ref{sec:Galois-CommRings}. Then the natural morphism of
Eilenberg-Mac~Lane spectra $HR\lra HS$ makes $HS/HR$ into a $G$-Galois
extension.
\end{ex}
\begin{ex}\label{ex:GalExtnMap=Iso}
Let $EG$ be any contractible space on which~$G$ acts freely. Let
$B/A$ be a $G$-Galois extension. Then $F(EG_+,B)/A$ is a $G$-Galois
extension and the collapse map $EG_+\lra S^0$ induces a morphism of
$A$-algebras $B\lra F(EG_+,B)$ commuting with the actions of~$G$ and
which is an equivalence of $B$-algebras.
\end{ex}

\begin{ex}\label{ex:KU/KO}
Let $\iota\:KO\lra KU$ be the complexification morphism which can be
given the structure of a morphism of commutative $\mathbb{S}$-algebras.

The action of $C_2=\<\gamma_2\>$ originates in the action of the stable
operation $\psi^{-1}$ whose action on $KU_{2n}=\Z u^n$ satisfies
\[
\gamma_2\.u^n=\psi^{-1}(u^n)=(-1)^nu^n.
\]
Recall that
\begin{equation}\label{eqn:KO*}
KO_*=\Z[\eta,y,w,w^{-1}]/(2\eta,\eta^3,y\eta,y^2-4w),
\end{equation}
where $\eta\in KO_1$, $y\in KO_4$ and $w\in KO_8$.

By the `Theorem of Reg Wood'~\cite[p.~206]{JFA-Chicago},
multiplication by the non-zero element $\eta\in KO_1$ induces a
cofibre sequence of $KO$-modules
\[
\Sigma KO\xrightarrow{\eta}KO\lra KO\wedge\cone(\eta)
\]
in which $KO\wedge\cone(\eta)\sim KU$ as $KO$-modules. This makes it
clear that $KU$ is self-dual in $\mathscr{D}_{KO}$ since in
$\mathscr{D}_{\mathbb{S}}$ we have $ D
\cone(\eta)\sim\Sigma^{-2}\cone(\eta)$. Then as $KU$-modules,
$KU\Smash_{KO}KU\sim KU\wedge\cone(\eta)$. Using this equivalence,
Rognes shows in~\cite{R:Opusmagnus} that (BNG-2) holds.

Notice that although $KU_*$ is not a projective module over $KO_*$,
$KU$ is in fact faithful over $KO$ since if $M$ is a $KO$-module with
$KU\Smash_{KO}M\sim*$, then using the cofibre sequence
\[
\Sigma KO\Smash_{KO}M\xrightarrow{\eta\wedge\id}KO\Smash_{KO}M
\lra KU\Smash_{KO}M
\]
we see that $\Sigma
KO\Smash_{KO}M\xrightarrow{\eta\wedge\id}KO\Smash_{KO}M\iso M$ is a
weak equivalence. Now since $\eta\in KO_1$ is nilpotent, this
implies that $M\sim*$.
\end{ex}

In the following, we state and prove some basic results about Galois
extensions which we will need later. These were stated by J.~Rognes
around~2000 and proofs can now be found in \cite[lemmas 6.1.2,
6.4.3, proposition 6.4.7]{R:Opusmagnus}.
\begin{thm}\label{thm:BraveNewGalois-Properties}
Let $A$ and $B$ be commutative $\mathbb{S}$-algebras and let $A\lra B$
be a $G$-Galois extension. Then in the derived category of $A$-modules
$\mathscr{D}_A$, the following hold.
\begin{enumerate}
\item[(a)]
$B$ is strongly self-dual.
\item[(b)]
For every $B$-module $N$, $N \Smash_A B \sim F(G_+,N)$.
\item[(c)]
For every $B$-module $N$, $N\wedge G_+\sim F_A(B,N)$.
In particular, $B\wedge G_+\sim F_A(B,B)$.
\end{enumerate}
\end{thm}
\begin{proof}
In the proof we will make extensive use of notions developed in
Section~\ref{sec:Duality}. The key part is~(a), the others follow
by formal arguments involving strong duality and we omit these.

The idea is to emulate as far as possible the ideas used in proving
the algebraic results of Section~\ref{sec:Galois-CommRings}. The
most important ingredient is a (weak) \emph{trace morphism}
\[
B\lra B^{\mathrm{h}G}\xrightarrow{\;\sim\;}A
\]
which factorizes the symmetrization map $\sum_{\gamma\in
G}\gamma\:B\lra B$. A construction for such a map can be found in
\cite[theorem~5.10]{JG&JM:Tate}. In our context (with $B_G$ denoting
the naive $G$-spectrum associated with~$B$), this produces a
homotopy commutative diagram
\begin{equation}\label{eqn:Tate-Diagram}
\xymatrix{
(B_G\wedge G_+)/G\iso B
\ar[rr]^{\sum_{\gamma\in G}\gamma}
\ar[d]_{(\id\wedge\inc)/G}
\ar@{-->}[drr]^{\tr}
&& B\sim F(EG_+,B) \\
(B_G\wedge EG_+)/G
\ar[rr]_{\bar\tau}
&&
A\sim B^{\mathrm{h}G}=F(EG_+,B_G)^G
\ar[u]_{\inc}
}
\end{equation}
and we take for our trace map the composition
\begin{equation}\label{eqn:Trace-Defn}
\tr=\tr_{B/A}=\bar\tau\o(\id\wedge\inc)/G\:B\lra A.
\end{equation}
It is straightforward to check that this is in fact a morphism
of $A$-modules. Having obtained a trace map, we can now define
an evaluation map to be the \emph{trace pairing}
\[
\epsilon\:B\Smash_AB\xrightarrow{\mult}B\xrightarrow{\tr}A
\]
which is a morphism of $A$-modules.
\begin{rems}\label{rem:tracemap}{\ }
\noindent
(i) If $B_G$ is $G$-equivariantly of the form $A\wedge G_+$, then
it is known from~\cite[proposition 2.4]{JG&JM:Tate} that the map
$\bar\tau$ in~\eqref{eqn:Tate-Diagram} is an isomorphism. \\
(ii) If $N\ideal G$, then there is a homotopy factorization
\[
\tr_{B/A}\sim\tr_{B/B^{\mathrm{h}N}}\tr_{B^{\mathrm{h}N}/A}.
\]
We are grateful to J.~Greenlees for showing us a verification of
this formula. \\
(iii) We claim that  the trace map $\tr\:B\lra A$ is $G$-invariant
in the sense that for any $\gamma\in G$, the composition
\[
B\xrightarrow{\gamma}B\xrightarrow{\tr}A
\]
is homotopic to $\tr$.

Using the description in~\cite[pp.~38--42]{JG&JM:Tate} we can
write the trace map as the composition
\[
\xymatrix{
{B \cong (B_G \wedge G_+)/G} \ar@{.>}[d]_{\tr}
\ar[rr]^{(\mathrm{id} \wedge i)/G} & &
{(B_G \wedge EG_+)/G} \ar[r]^{\ph{\sim}\widetilde{\tau}\ph{\sim}}
& {i^*i_*(B_G \wedge EG_+)^G} \ar[d]^{\varepsilon}
\\
{F(EG_+, B_G)^G} \ar[rr]^{\sim} & & {F(EG_+,i^*i_*(B_G))^G}
& {i^*i_*(B_G)^G} \ar[l]_{\ph{\sim\sim\sim}\varepsilon^G}
}
\]
We should mention that there is a $G$-action on $(B_G \wedge EG_+)/G$.
We consider the semidirect product $G \ltimes G$ where we take the
conjugation action of $G$ on itself. Then $G \ltimes G$ acts on
$B_G \wedge EG_+$, and if we divide out by the normal subgroup,
we are still left with a $G$-action on the quotient. This is the
$G$-action that Greenlees and May use on $(B_G \wedge EG_+)/G$
(see~\cite[p.~38]{JG&JM:Tate}).

The map $\widetilde{\tau}$ is a transfer map and is natural, hence
it is equivariant. We denote by $\varepsilon$ maps induced by the
collapse map $EG_+ \lra S^0$. The last map is the unit of the
adjunction $(i^*,i_*)$ defined in~\cite[lemma~0.1]{JG&JM:Tate},
this is an equivariant map although it is not a weak equivalence
in the equivariant setting. Therefore the trace is $G$-invariant
as claimed and hence the self-duality of
Theorem~\ref{thm:BraveNewGalois-Properties}(a) is given by a
$G$-equivariant equivalence $B \lra F_A(B,A)$. This map is adjoint
to the trace pairing
$$
B \wedge_A B \xrightarrow{\mult}  B \xrightarrow{\tr} A
$$
which is clearly equivariant if we take the diagonal $G$-action
on $B \wedge_A B$.
\end{rems}

We also need to produce a coevaluation $\eta\:A \lra B\Smash_AB$.
Working in the derived category $\mathscr{D}_A$, this is done
using the map $B\lra B\Smash_AB$ implicit in Condition~(BNG-2)
of Definition~\ref{defn:BraveNewGaloisExtn-Rognes} and splitting
the multiplication map $B\Smash_AB\lra B$ which corresponds
to projection onto the identity element component of
$\prod_{\gamma\in G}B\iso F(G_+,B)$. The composition
$\eta\:A\lra B\lra B\Smash_AB$ can be viewed as the unique
element of $\pi_0(B\Smash_AB)$ projecting to the element
of $(\delta_{\gamma,1})\in\prod_{\gamma\in G}\pi_0B$, where
\[
\delta_{\alpha,\beta}=
\begin{cases}
1& \text{if $\alpha=\beta$}, \\
0& \text{otherwise}.
\end{cases}
\]

Now to show that $B$ is strongly self-dual with evaluation $\epsilon$
and coevaluation $\eta$, we have to verify Condition~(b) of
Theorem~\ref{thm:Duality}. We need to check that when $P=B=Q$,
the compositions in the diagrams~\eqref{eqn:Duality} are indeed
the identity morphisms.

%

For each $\gamma\in G$, there is a map of $A$-ring spectra
\begin{equation}\label{eqn:EtaleMap-factor}
\tilde\gamma\:B\Smash_AB\xrightarrow{\id\wedge\gamma}B\Smash_AB
\xrightarrow{\mu}B.
\end{equation}
Recalling that as maps from $B$ to $B$,
\[
\iota\o\tr\sim\sum_{\gamma\in G}\gamma,
\]
we find that the composition in~\eqref{eqn:Duality-P} is
\[
\sum_{\gamma\in G}\tilde\gamma\o(\id\wedge\mult)\o(\eta\wedge\id)
\:B\iso A\Smash_AB\lra B.
\]
Since each $\tilde\gamma$ is a $B$-bimodule morphism and because
of the way $\eta$ was characterized in terms of its projections
under the $\tilde\gamma$, this composition is homotopic to
\[
\sum_{\gamma\in G}\delta_{\gamma,1}\gamma=\id.
\]
A similar discussion applies to the composition in~\eqref{eqn:Duality-Q}.
\end{proof}

\begin{rem}\label{rem:tr-NotEpi}
In general, in distinction to part~(b) of
Theorem~\ref{thm:GaloisExtn-Properties}, the trace map $\tr=\tr_{B/A}$
need not induce an epimorphism $\tr_*\:B_*\lra A_*$. For example, in the
case of $KU/KO$ discussed in Example~\ref{ex:KU/KO},
the trace map $\tr$ agrees with the realification map and
$\tr_*(u)=\eta^2$. To see this, note that we are dealing with
$KO$-module maps $KU\lra KO$. Then by the self-duality of the
$KO$-module $KU$,
\begin{equation}\label{eqn:KU/KOtr=realization}
\pi_0F_{KO}(KU,KO)\iso\pi_0F_{KO}(KO,KU)\iso\pi_0KU\iso\Z.
\end{equation}
This means that elements of $\mathscr{D}_{KO}(KU,KO)\iso\pi_0F_{KO}(KU,KO)$
are detected by their induced effect on the free abelian homotopy groups
$\pi_{4n}KU\lra\pi_{4n}KO$ for $n\in\Z$, and using the $KO_*$-module structure
we find that this is determined by the homomorphism $\pi_0KU\lra\pi_0KO$.
Thus it suffices to know the standard fact that the realification map
$KU\lra KO$ induces $2\:\pi_0KU\lra\pi_0KO$.
\end{rem}

When the order of the Galois group $G$ is invertible in $A_*$, such
anomalies do not occur.
\begin{prop}\label{prop:BraveNewGaloisExtn-easycase}
Let $A$ and $B$ be commutative $\mathbb{S}$-algebras and let $A\lra B$
be a weak $G$-Galois extension for which $A_*$ is a $\Z[1/|G|]$-algebra.
Then the following hold.
\begin{enumerate}
\item[(a)]
The unit $A\lra B$ induces a monomorphism $A_*\lra B_*$.
\item[(b)]
The unit $A \lra B$ and trace $\tr\:B\lra A$ compose to an equivalence
on $A$.
Hence an $A$-module $M$ is a retract of $B\Smash_AM$. In particular,
$B$ is faithful.
\end{enumerate}
\end{prop}
\begin{proof}
(a) There is a spectral sequence
\begin{equation}\label{eqn:FixedPtSS}
\mathrm{E}_2^{s,t}=\mathrm{H}^s(G;B_t)\Lra(B^{\mathrm{h}G})_{t-s}.
\end{equation}
Under the above hypotheses, the $\mathrm{E}_2$-term is concentrated
in the zero-line where $\mathrm{E}_2^{0,t}= (B_t)^G$.
Hence on homotopy groups, the unit induces the inclusion of the fixed
points of $B_*$. \\
(b) In the diagram of~\eqref{eqn:Tate-Diagram}, precomposition of the
top row with the unit $\iota\:A\lra B$ induces an equivalence
$\tr\iota\:A\lra A$. Thus the trace and unit are split. It follows that
an $A$-module $M$ is a retract of $B\Smash_AM$.
\end{proof}

Part~(b) of the last Proposition implies the following.
\begin{cor}\label{cor:BraveNewGaloisExtn-easycase}
Let $A$ and $B$ be commutative $\mathbb{S}$-algebras and let $B/A$ be
a weak $G$-Galois extension for which $A_*$ is a $\Z[1/|G|]$-algebra.
Then $B/A$ is a $G$-Galois extension.
\end{cor}

The following example is  straightforward to verify. Let~$p$ be an
odd prime. The Johnson-Wilson spectrum $E(1)$ agrees with the Adams
summand of $KU_{(p)}$. Passing to the $p$-completions, there is an
action of $C_{p-1}$ on $(KU_p)$ (see~\cite[theorem 9.2]{AB&BR},
\cite[\S 7]{PG&MH}) that is easily seen to turn $KU_p$ into a weak
Galois extension of $E(1)_p$. We can use
Corollary~\ref{cor:BraveNewGaloisExtn-easycase} to obtain
\begin{ex}\label{ex:KU-G}
$KU_p/E(1)_p$ is a $C_{p-1}$-Galois extension.
\end{ex}

There is analogue of Proposition~\ref{prop:GAL-Groupoid}, namely
\begin{prop}\label{prop:GalExtnMap=Iso}
Let $B/A$ and $C/A$ be $G$-Galois extensions and let $\phi\:B\lra C$
be a morphism of $A$-algebras which commutes with the action of\/ $G$.
Then $\phi$ is a weak equivalence.
\end{prop}
\begin{proof}
By the faithfulness of $B$ and $C$, it suffices to check
this for the morphism $\tilde\phi=\id\wedge\id\wedge\phi$
between the $G$-Galois extensions $B\Smash_AC\Smash_AB/B\Smash_AC$
and $B\Smash_AC\Smash_AC/B\Smash_AC$. But then
\[
B\Smash_AC\Smash_AB \sim \prod_{\gamma\in G}B\Smash_AC
\sim B\Smash_AC\Smash_AC.
\]
Now by construction, the map $\id \wedge \Theta$ of
Definition~\ref{defn:BraveNewGaloisExtn-Rognes}(BNG-2)
is $B \wedge C$-linear, hence
$\tilde \phi$ induces a $B^A_*C$-algebra homomorphism
\[
\tilde\phi_*\:\prod_{\gamma\in G}B^A_*C\lra\prod_{\gamma\in G}B^A_*C
\]
which is $G$-equivariant and by Proposition~\ref{prop:GAL-Groupoid},
this is an isomorphism.
\end{proof}

\part{From algebraic Galois extensions to brave new Galois extensions}
\label{part:2}

\section{Topological realization of algebraic Galois extensions}
\label{sec:TopRealGalExtnRngs}

Let $A$ be a commutative $\mathbb{S}$-algebra  and let $G$ be
a finite group. Also recall Proposition~\ref{prop:GAL-Groupoid}.
\begin{thm}\label{thm:A->B}
If $B_*$ is a\/ $G$-Galois extension of $A_*$, then there is
a commutative $A$-ring spectrum $B$ realizing $B_*$ as $\pi_*B$
and a homotopy action of\/ $G$ on $B$ by morphisms of $A$-ring
spectra which induce the action of\/ $G$ on $B_*$. Furthermore,
if $C_*$ is also a $G$-Galois extension of $A_*$ and there is
a $G$-isomorphism $\Phi\:B_*\lra C_*$ of $A_*$-algebras, then
there is a map of $A$-ring spectra $\phi\:B\lra C$ which induces
$\Phi$. It is $G$-equivariant up to homotopy.
\end{thm}
\begin{proof}
By Theorem~\ref{thm:GaloisExtn-Eqces}(d), $B_*$ is a finitely
generated projective $A_*$-module, so we can realize $B_*$ as
the image of an idempotent
$e\: \bigoplus_{i=1}^n \Sigma^{m_i} A_*
                      \lra \bigoplus_{i=1}^n \Sigma^{m_i} A_*$.
We can model the map $e$ on a wedge of suspensions of $A$.
Therefore the mapping telescope of
\begin{equation}\label{eqn:adamscond}
\bigvee \Sigma^{m_i} A \xrightarrow{e} \bigvee \Sigma^{m_i} A
\xrightarrow{e} \cdots
\end{equation}
gives rise to an $A$-module spectrum $B$ with $\pi_*B \cong B_*$.
%

The K\"unneth spectral sequence
\begin{equation}\label{eqn:KSS}
\mathrm{E}^2_{p,q}=\Tor^{A_*}_{p,q}(B_*,B_*)\Lra
B^{\ph{p}A\ph{q}}_{p+q}B
\end{equation}
of~\cite{EKMM} collapses to give
\begin{equation}\label{eqn:KSS-BsmashB}
B^A_*B=B_*\oTimes_{A_*}B_*.
\end{equation}
More generally, for each $n\geq2$,
\begin{equation}\label{eqn:KSS-B^n}
\pi_*B^{(n)} \iso \overset{n-1}
{\overbrace{B^A_*B\oTimes_{B_*}B^A_*B\oTimes_{B_*}\cdots\oTimes_{B_*}B^A_*B}},
\end{equation}
where $(\ph{B})^{(n)}$ denotes the $n$-fold smash product over~$A$.
This is projective both as an $A_*$-module and as a $B_*$-module.

{}From~\cite{EKMM}, for each $A$-module $Y$ there is a universal
coefficient spectral sequence
\begin{equation}\label{eqn:UCSS-BProduct}
\mathrm{E}_2^{p,q}= \Ext_{A_*}^{p,q} (\pi_*B^{(n)},Y_*) \Lra
Y_A^{p+q}(B^{(n)}).
\end{equation}
By the projectivity of the first variable, this spectral sequence
collapses to give
\begin{equation}\label{eqn:B-Product}
Y_A^*(B^{(n)}) \iso\Hom_{A_*}^* (\pi_*B^{(n)},Y_*).
\end{equation}
The product on $B_*$ is an element of
$\Hom_{A_*}(\pi_*(B\Smash_A B),B_*)$ which corresponds
to a unique element of $B^0(B\Smash_A B)$. Since $B_*$
is a commutative $A_*$-algebra, this product on $B$ is
homotopy associative, commutative and unital over~$A$.

Similarly, the action of elements of $G$ on $B_*$ induces
a homotopy action of~$G$ on~$B$ by morphisms of $A$-ring
spectra.

Since $B_*$ is finitely generated and projective as an
$A_*$-module, the relevant universal coefficient spectral
sequence collapses and for $m,k\in\N$, $n\in\Z$ there are
isomorphisms
\[
\mathscr{D}_A(B^{(k)},\Sigma^nC^{(m)})\xrightarrow{\iso}
\Hom_{A_*}^n((B_*)^{\otimes k},(C_*)^{\otimes m}),
\]
where $(\ph{B})^{(k)}$ denotes the $k$-fold smash product over~$A$.

Let $\Hom_{A_*\mathrm{-alg}}(B_*, C_*)$ denote the set
of $G$-equivariant $A_*$-algebra maps from $B_*$ to $C_*$,
and let
$\Hom_{A_*\mathrm{-alg}}(B_*, C_*)^G\subset\Hom_{A_*\mathrm{-alg}}(B_*, C_*)$
be the subset of $G$-equivariant maps. The latter can
be written as an iterated equalizer: first we obtain
$\Hom_{A_*\mathrm{-alg}}(B_*, C_*)$ as the equalizer
$$
\xymatrix{ {\Hom_{A_*\mathrm{-alg}}(B_*, C_*)} \ar[r] &
{\Hom_{A_*}(B_*, C_*)} \ar@<0.5ex>[r] \ar@<-0.5ex>[r] &
{\Hom_{A_*}(B_* \otimes_{A_*} B_*, C_*)}
}
$$
The projectivity of $B_*$ over $A_*$ therefore gives
that $\Hom_{A_*\mathrm{-alg}}(B_*, C_*)$ is the same
as the homotopy classes of maps of $A$-ring spectra
from $B$ to $C$, $[B,C]_{A\mathrm{-ring}}$. Similarly,
as $G$ is finite we know that $(G_+ \wedge B)_*$ is
finitely generated projective over $A_*$ and we obtain
that the homotopy classes of $G$-equivariant maps of
$A$-ring spectra, $[B,C]_{A\mathrm{-ring}}^G$, are the
same as $\Hom_{A_*\mathrm{-alg}}(B_*, C_*)^G$. The map
$\Phi$ is an element of the latter; therefore there is
a realization $\phi\:B \lra C$, which is homotopy
$G$-equivariant and a map of $A$-ring spectra.
\end{proof}

\section{Realizing Galois extensions}\label{sec:RealisingGalExtns}

In this section we will assume that the following two conditions
hold.
\begin{enumerate}
\item[(GE-1)]
$A$ is a commutative $\mathbb{S}$-algebra and $B$ is a commutative
$A$-ring spectrum. There is a homotopy action of the finite
group~$G$ on $B$ viewed as an $A$-module, \ie, there is a
homomorphism of monoids $G\lra\mathscr{D}_A(B,B)$ which is in fact
an action by automorphisms of~$B$ as an $A$-ring spectrum, \ie, the
homomorphism $G\lra\mathscr{D}_A(B,B)$ is compatible with the
product $\mu\:B\Smash_A B\lra B$ in $\mathscr{D}_A$.
\item[(GE-2)]
The $A_*$-algebra $B_*$ is a $G$-Galois extension with respect
to the induced action of~$G$, thus the action of~$G$ on $B_*$ is
effective, \ie, $G \leq \Aut_{A_*}(B_*)$,  and satisfies the axioms
of a Galois action (see Definition~\ref{defn:GaloisExtn}).
\end{enumerate}
For instance, these conditions are satisfied when we start with
a situation as in Theorem~\ref{thm:A->B}.

\begin{prop}\label{prop:BSmashB}
As an $A_*$-algebra,
\[
B^A_*B=\pi_*(B\Smash_AB) \cong \prod_{\gamma\in G}B_*,
\]
where the map is induced by that of~\eqref{eqn:EtaleMap}. Hence
$B^A_*B$ is an \'etale $B_*$-algebra.
\end{prop}
\begin{proof}
Recalling~\eqref{eqn:KSS-BsmashB}, we see that
\[
B^A_*B\iso B_*\oTimes_{A_*}B_*.
\]
As the edge homomorphism in the K\"unneth spectral sequence
\eqref{eqn:KSS} is multiplicative, this is an isomorphism of
$A_*$-algebras. Since $B_*$ is \'etale over $A_*$,
$B_*\oTimes_{A_*}B_*$ is \'etale over $B_*$. Hence $B^A_*B$
is \'etale over $B_*$.
\end{proof}
\begin{cor}\label{cor:BSmashB}
For any $B^A_*B$-bimodule $M_*$ and $B^A_*B$-module $N_*$, the
Hochschild cohomology and the $\Gamma$-cohomology of $B^A_*B$
relative to $B_*$ vanish, \ie,
\[
\HH^{**}(B^A_*B\mid B_*;M_*) \cong  M_*
\quad\text{ and }\quad
\HG^{**}(B^A_*B\mid B_*;N_*) = 0.
\]
\end{cor}
\begin{prop}\label{prop:Obstructions-a}
Assume that $A,B,G$ satisfy conditions {\rm(GE-1)} and {\rm(GE-2)}.
Then the $A$-ring spectrum structure on $B$ has a refinement to a
commutative $A$-algebra structure which is unique up to contractible
choice.
\end{prop}
\begin{proof}
We will use adaptations of the obstruction theory of Robinson to
the relative case. This approach has been set up to establish the
existence of $E_\infty$ structures (or equivalently commutative
$\mathbb{S}$-algebra structures) on spectra. Our aim is to establish
commutative $A$-algebra structures on a homotopy commutative
$A$-ring spectrum~$B$.

The geometric nature of Robinson's  obstruction groups as described
in~\cite[\S5]{Ro:Einfty} ensures that the obstructions for imposing
a commutative $A$-algebra structure on~$B$ live in
$\Gamma$-cohomology $\HG^{**}(B^A_*B\mid B_*; B_*)$ of
$(B\Smash_A B)_*$ relative to $B_*$. Using the notation
of~\cite[definition~5.3]{Ro:Einfty}, an $n$-stage for such a structure
corresponds to action maps
\[
\mu_m\:\nabla^n\mathcal{T}_m \ltimes_{\Sigma_m}
\overset{m}{\overbrace{B\Smash_A\cdots\Smash_AB}}\lra B
\]
for $m \leq n$ and certain compositions. As we have assumed that
$A_*\lra B_*$ is $G$-Galois, we find that $B_*$ is $A_*$-projective;
hence one obtains universal coefficient and K\"unneth isomorphisms
which identifies $B^*_A(B\wedge_A B)$ with
$\Hom_{B_*}(B_*\otimes_{A_*}B_*,B_*)$.

The result is now immediate using the modified obstruction theory
of Robinson~\cite{Ro:Einfty,Ro:unit} and general properties of
$\Gamma$-cohomology established by Robinson and
Whitehouse~\cite{RoWh} for $A\lra B$, because $B_*\otimes_{A_*}B_*$
is \'etale over $B_*$ and hence the obstruction groups vanish.
\end{proof}

Now by making use of results of Robinson and
Whitehouse~\cite{Ro:Einfty,Ro:unit,RoWh}, and Goerss and
Hopkins~\cite{PG&MH}, we obtain
\begin{thm}\label{thm:Obstructions}
Assume that $A,B,G$ satisfy conditions {\rm(GE-1)} and {\rm(GE-2)}.
Then the following hold.
\begin{enumerate}
\item[(a)]
Each element $\gamma\in G$ induces a morphism of $A$-algebras from
$B$ to $B$ which is unique up to contractible choice.
\item[(b)]
The morphisms of part {\rm(a)} combine to give an action of\/ $G$
on~$B$ by $A$-algebra automorphisms.
\item[(c)]
Suppose that $A,C,G$ also satisfy the conditions of {\rm(GE-1)} and
{\rm(GE-2)}, thus there is a unique $A$-algebra structure on $C$ as
in {\rm(a)}. If $\tilde{\phi}\: B \lra C$ is a map of $A$-ring
spectra which is $G$-equivariant up to homotopy, then there are
commutative $A$-algebras $B''$ and $C''$ together with weak
equivalences $B \sim B''$ and $C \sim C''$. These weak equivalences
are zigzags of weak equivalences of commutative $A$-algebras which
are $G$-equivariant up to homotopy. There is a map of commutative
$A$-algebras $\bar{\phi}: B'' \ra C''$ which is strictly
$G$-equivariant and which induces $\tilde{\phi}$.
\end{enumerate}
\end{thm}
\begin{proof}
For (a) \& (b), the desired result comes from the triviality of the
spectral sequence for the homotopy of the derived space of
commutative $A$-algebra self-maps of~$B$ based at any choice of map
in the $A$-algebra homotopy class of an element of~$G$ realized as
an $A$-algebra morphism. More precisely, we use the generalization
of~\cite[theorem~4.5]{PG&MH} to the setting of $A$-algebras and take
$E=X=Y=B$. Note that the description of $B$ as in
\eqref{eqn:adamscond} ensures that $B$ satisfies the Adams condition
required in \cite[definition 3.1]{PG&MH}.

Now we make a modification of~\cite[definition 3.2]{PG&MH}, using as
$\mathcal{P}$ the set of spectra consisting of the $A$-sphere $S_A$,
$B$, their suspensions and finite wedges of these. As $B_*$ is
finitely generated $A_*$-projective, we have a universal coefficient
isomorphism.

The second quadrant spectral sequence converging to the homotopy
groups of the derived space of self-maps of $B$ in the category of
$E_\infty$-algebras in $A$-modules,
$\map_{E_\infty\mbox{-}A\mathrm{\mbox{-}alg}}(B,B)$, looks as follows.
The $\mathrm{E}_2$-term is
\[
\mathrm{E}_2^{s,t} =
\begin{cases}
\Hom_{B_*\mathrm{\mbox{-}alg}}(B^A_*B,B_*) & \text{for $(s,t)=(0,0)$}, \\
\Der^s_{B_*}(B^A_*B, \Omega^tB_*) & \text{for $t > 0$},
\end{cases}
\]
where $\Der^s_{B_*}(B^A_*B,\Omega^tB_*)$ denotes the $s$-th derived
functor of derivations into the $t$-th shift of $B_*$ and
$\Hom_{B_*\mathrm{\mbox{-}alg}}(B^A_*B,B_*)$ are the morphisms of
$B_*$-algebras from $B^A_*B\cong B_*\otimes_{A_*}B_*$ to $B_*$.
For $s>0$ we know that $\Der^s_{B_*}(B^A_*B,\Omega^tB_*)$ vanishes
since $B^A_*B$ is \'etale over $B_*$. The reader might wish to use
the comparison result of~\cite{MB&BR} to see that. In this case it
provides an isomorphism between  $\Der^s_{B_*}(B^A_*B, \Omega^{-t}B_*)$
and $\HG^{s,t}(B_*^AB|B_*;B_*)$. There should be a direct proof as
well. As a basepoint for the derived mapping space
$\map_{E_\infty\mbox{-}A\mathrm{\mbox{-}alg}}(B,B)$ we take the
identity map of~$B$.

Since $B_*$ is $G$-Galois over $A_*$, each group element $\gamma\in G$
gives rise to an element in the morphism set
$\Hom_{B_*\mathrm{\mbox{-}alg}}(B^A_*B,B_*)$ which sends $b_1\otimes b_2$
to $b_1{}^\gamma b_2$.

Using the partially invisible corollary~\cite[Corollary to theorem~4.5,
analogous to~4.4]{PG&MH}, we see that the vanishing of the obstruction
groups $\Der^s_{B_*}(B_*^AB,\Omega^sB_*)$ and
$\Der^{s+1}_{B_*}(B_*^AB,\Omega^sB_*)$ implies that the Hurewicz map
\[
\pi_0(\map_{E_\infty\mbox{-}A\mathrm{\mbox{-}alg}}(B,B))
\lra
\Hom_{B_*\mathrm{\mbox{-}alg}}(B_*^AB, B_*)
\]
is a bijection, hence each component is labelled by an element of
the latter. (The condition for surjectivity is stated
in~\cite[theorem~4.5]{PG&MH}.) Therefore the derived mapping space
$\map_{E_\infty\mbox{-}A\mathrm{\mbox{-}alg}}(B,B)$ has contractible
components, because for  \'etale algebras like $B^A_*B$ Gamma
cohomology vanishes with arbitrary coefficients. In addition we see
that the group $G$ is a submonoid of
$\pi_0(\map_{E_\infty\mbox{-}A\mathrm{\mbox{-}alg}}(B,B))$, in
particular every $\gamma \in G$ gives a self-map of $B$ in the
homotopy category of $E_\infty$-$A$-algebras. In the terminology
of~\cite{DKS} (or of~\cite[definition 2.1]{DH}) the diagram category
consisting of the group viewed as a one-object category gives rise
to an $h_\infty$-diagram: Let $X(0)$ be $B$ and for every $\gamma\in G$
we obtain a self-map $X(0,\gamma)$ of $B$. The path-component of the
image of $\gamma$ is contractible. Using~\cite[theorem~2.2]{DH} we
can strictify this diagram such that there is a weakly equivalent
$E_\infty$-$A$-algebra $B'$ with a strict $G$-action.

For part (c) the arguments are a little bit more involved. As we
saw, we can realize $A$-algebras $B'$ and $C'$ with  actual
$G$-actions, such that $B'$ is weakly equivalent to $B$ and $C'$
is weakly equivalent to~$C$ via maps of $E_\infty$-$A$-algebras.
As $\tilde{\phi}$ was a map of $A$-ring spectra, it gives rise
to a map $\phi$ of homotopy $A$-ring spectra from $B'$ to $C'$
which is $G$-equivariant up to homotopy. In particular, the map
on $C$-homology, $C_*(\phi)$ is a map of commutative $C_*$-algebras.
As $B_*$ is \'etale over  $A_*$ a base-change argument implies
that $C_*^AB \cong C_* \otimes_{A_*} B_*$ is $C_*$-\'etale.

Using the Hurewicz argument again, $C_*(\phi)$ gives rise to a
map $\phi'$ of $E_\infty$-$A$-algebras from  $B'$ to $C'$ which
is still $G$-equivariant up to homotopy.

We claim that the following is an $h_\infty$-diagram: take
$X(0) = B'$ and $X(1)=C'$ as vertices and for every group element
$\gamma\in G$ we get morphisms $X(\gamma,0)$ from $X(0)$ to itself
and $X(\gamma,1)$ on $X(1)$. We place the map $\phi'\:B' \lra C'$
in the diagram as a connection between $X(0)$ and $X(1)$:
$$
\xymatrix{{}&{}&{}&{}&{}\\
{}&{} &{X(0)=B'}\ar[r]^{\phi'} \ar@(ur,ul)_{X(\gamma,0)}
\ar@(dl,dr)_{X(\gamma',0)} \ar`l[u]`u[r]^{X(\gamma'',0)}`r[dl]
& {C'=X(1)} \ar@(ur,ul)_{X(\gamma,1)} \ar@(dl,dr)_{X(\gamma',1)}
\ar@{-<}`r[u]`[l]_{X(\gamma'',1)}`[dr]
& {}
\\
{}&{}&{}&{}{}& }
$$
For the element in
$\pi_0(\map_{E_\infty\mbox{-}A\mathrm{\mbox{-}alg}}(B,C))$
corresponding to $C_*(\phi)$ we get an arrow from $X(0)$ to $X(1)$
in the homotopy category of $E_\infty$-$A$-algebras. As every
component in $\map_{E_\infty\mbox{-}A\mathrm{\mbox{-}alg}}(B,C)$
is contractible, one can strictify $X$ to get a weakly equivalent
diagram $\bar{X}$. By construction, these equivalences are maps of
commutative $A$-algebras which are $G$-equivariant up to homotopy.
The arrow $\bar{\phi}$ from $\bar{X}(0)$ to $\bar{X}(1)$ is
$G$-equivariant by construction and is a weak equivalence due to
three-out-of-four.
\end{proof}
%
%
\begin{prop}\label{prop:NormalBasis-Trace}
Assume that $A,B,G$ satisfy conditions {\rm(GE-1)} and {\rm(GE-2)}.
Then the trace map $\tr:B\lra A$ induces a surjection
$\tr_*\:B_*\lra A_*$.
\end{prop}
\begin{proof}
As $B_*/A_*$ is $G$-Galois, the trace map is an epimorphism by
Theorem~\ref{thm:GaloisExtn-Properties}(b).
\end{proof}
\begin{cor}\label{cor:NormalBasis-Trace}
The unit $A\lra B$ is split. Hence, $B$ is faithful as an
$A$-module.
\end{cor}
\begin{proof}
Let $b\:A\lra B$ be an element of $B_0=\pi_0B$ for which
$\tr_*(b)=1$. Then the composition
\[
B\sim A\Smash_AB\xrightarrow{b\wedge\id}B\Smash_AB
                   \xrightarrow{\mult}B\xrightarrow{\tr}A
\]
splits the unit. Hence for any $A$-module $M$, $M$ is a retract
of $B\Smash_AM$ and so $B$ is faithful.
\end{proof}

Now we can state and prove the main result of this section.
\begin{thm}\label{thm:AB-New}
Assume that $A,B,G$ satisfy conditions {\rm(GE-1)}, {\rm(GE-2)}.
Then $B/A$ is a $G$-Galois extension.
\end{thm}
\begin{proof}
Without loss of generality we can replace $B$ by a cofibrant
commutative $A$-algebra: there is a functorial cofibrant replacement
functor $Q(-)$ (see \cite[VII, \S \S 4,5]{EKMM}); therefore $Q(B)$
inherits the $G$-action from $B$. As $\pi_*Q(B) \cong B_*$ is still
$A_*$-projective, $Q(B)$ is unramified. In the following we write
$B$ instead of $Q(B)$.

We use the homotopy fixed point spectral sequence
\[
\mathrm{E}_2^{s,t}=\mathrm{H}^s(G;B_t)\Lra(B^{\mathrm{h}G})_{t-s}
\]
to ensure that $B$ has the correct homotopy fixed points with respect
to the $G$-action. We suppress the internal grading to ease notation.

Now by Proposition~\ref{prop:G-Cohom-S/R},
\[
\mathrm{E}^*_2=\mathrm{E}^0_2=\Hom_{A_*G}(A_*, B_*)\iso(B_*)^G=A_*.
\]
Therefore $\pi_*(B^{\mathrm{h}G})\cong A_*$ and so
$B^{\mathrm{h}G}\sim A$.
\end{proof}
\begin{ex}\label{ex:Cyclotomic}
For $n\geq1$ and $\zeta_{n^\ell}$ a primitive $n^\ell$-th root of
unity, we may consider the ring $\Z[1/n,\zeta_{n^\ell}]\subset\C$.
The prime factors of the discriminant of $\Z[1/n,\zeta_{n^\ell}]$
over $\Z[1/n]$ are known to divide~$n$, so $\Z[1/n,\zeta_{n^\ell}]$
is unramified over $\Z[1/n]$. Then $\Z[1/n,\zeta_{n^\ell}]/\Z[1/n]$
is a $(\Z/{n^\ell})^\times$-Galois extension, and there is an
isomorphism of $\Z[1/n,\zeta_{n^\ell}]$-algebras
\[
\Z[1/n,\zeta_{n^\ell}]\otimes_{\Z[1/n]}\Z[1/n,\zeta_{n^\ell}]
      \iso \prod_{\gamma\in(\Z/{n^\ell})^\times}\Z[1/n,\zeta_{n^\ell}].
\]
For any commutative $\mathbb{S}$-algebra $A$, by Theorem~\ref{thm:A->B}
we can give $B=A\Z[1/n,\zeta_{n^\ell}]$ the structure of a commutative
$A$-ring spectrum. By applying Theorem~\ref{thm:Obstructions}, we see
that the ring $A_*\otimes\Z[1/n,\zeta_{n^\ell}]$ can be realized as
the homotopy ring of a commutative $A[1/n]$-algebra. Thus we find
that $B/A$ is a $(\Z/{n^\ell})^\times$-Galois extension. This gives
a different approach to results of~\cite{SVW}.
\end{ex}
\begin{ex}\label{ex:En}
For a prime $p$, let $E_n$ denote the $2$-periodic Lubin-Tate
spectrum whose homotopy ring is
\[
(E_n)_*= \mathrm{W}\mathbb{F}_{p^n}[[u_1,\ldots,u_{n-1}]][u,u^{-1}]
\]
where the $u_i$ are of degree zero and $u$ is an element of
degree~$-2$. This is known to be an algebra over the $I_n$-adic
completion $\hat{E(n)}$ of the Johnson-Wilson spectrum $E(n)$
(see~\cite{AB&BR} for a proof that $\hat{E(n)}$ is commutative).
On coefficients, the ring map from $\hat{E(n)}_*$ to $(E_n)_*$
is determined by $v_i \mapsto u_iu^{1-p^i}$. Then $E_n/\hat{E(n)}$
is a $C_n\ltimes\F_{p^n}^\times$-Galois extension.
\end{ex}
\section{Topological Harrison groups}\label{sec:TopHarGps}

As we want to compare algebraic Galois extensions to topological
ones, we propose the following definitions of  Harrison sets for
a commutative $\mathbb{S}$-algebra $A$ and a finite group $G$.
\begin{defn}\label{defn:weaktopHarreqt}
We call two $G$-Galois extensions of $A$,  $B'$ and $B''$,
\emph{weakly Harrison equivalent} if there are commutative
$A$-algebras with a homotopy $G$-action $(B_i)_{i=1}^n$ and
commutative $A$-algebras with strict $G$-action $(B'_i)_{i=1}^{n-1}$
with weak equivalences of $A$-algebras $\epsilon_i, \rho_i$ which
are homotopy $G$-equivariant as in the following diagram.
\begin{equation}\label{eqn:zigzag}
\xymatrix{
{B'}&{}&{B'_1}&{}&{}&{}&{B'_{n-1}}&{}&{B''}\\
{}&{B_1}\ar[ul]^{\epsilon_1}\ar[ur]_{\rho_1}
&{}&{B_2}\ar[ul]^{\epsilon_2}\ar[ur]_{\rho_2}&{}\ar@{}[ur]
  |{\cdots}&{} \ar[ur]_{\rho_{n-1}}&{}&{B_n}\ar[ul]^{\epsilon_n}\ar[ur]_{\rho_n}&{}
}
\end{equation}
We denote the set of such equivalence classes by $\Har^w(A,G)$.
\end{defn}
\begin{defn}\label{defn:topHarreqt}
We call two $G$-Galois extensions of $A$,  $B'$ and $B''$,
\emph{Harrison equivalent} if there are commutative $A$-algebras
with a strict $G$-action $(B_i)_{i=1}^n$ and commutative
$A$-algebras with strict $G$-action $(B'_i)_{i=1}^{n-1}$ with
weak equivalences of $A$-algebras $\epsilon_i, \rho_i$ which are
$G$-equivariant and which fit diagram \eqref{eqn:zigzag}. Such
equivalence classes are denoted by $\Har(A,G)$.
\end{defn}
\begin{rem}\label{rem:TopHarrison}
\begin{itemize}
\item[]
\item
The weak equivalence notion comes out of our realization result. We
propose the definition of $\Har(A,G)$ as a compromise. It is strong
enough to prove structural results. Note that there is an obvious map
$$ \Har(A,G) \lra \Har^w(A,G).$$
\item
A map in the other direction would require to replace $A$-algebras
with homotopy $G$-actions by $A$-algebras with strict $G$-action.
This means that we have to compare the mapping space of $G$-equivariant
maps from $EG$ to the mapping space of $A$-algebra endomorphisms
of an algebra spectrum $B$ with the mapping space of $G$-equivariant
maps from a point to that endomorphism space. This boils down to proving
the Sullivan conjecture for the endomorphism space of $A$-algebra maps
on $B$.
\end{itemize}
\end{rem}

\begin{thm}\label{thm:Harrison-Natural}
The constructions of the last section produce a map
\begin{equation}\label{eqn:Harrison-Natural}
\Real_{G}\:\Har(A_*,G)\lra \Har^w(A,G).
\end{equation}
\end{thm}
\begin{proof}
We must show that the topological realization of an algebraic
map $\Phi\: B_* \lra C_*$ of $G$-Galois extensions is unique
up to $G$-equivariant homotopy. Assume there are two such
realizations, $\phi$ and $\psi$ of $\Phi$, which are maps of
commutative $A$-algebras. The connected components of
$\map_{E_\infty\mbox{-}A\mathrm{\mbox{-}alg}}(B,C)$ are
labelled by the elements of
$\Hom_{A_*\mathrm{-alg}}(B_*,C_*)$. Thus $\phi$ and $\psi$ are
in the same path-component and a path between these gives a
homotopy $H$. The group action of $G$ preserves the components.
Therefore any element $g \in G$ applied to $\phi$ and $\psi$ is
again in that component and so is $g$ applied to $H$. We connect
$\phi$ with $g \phi g^{-1}$, and similarly $\psi$ and $H$ with
their conjugates. Schematically this yields the following diagram.
$$
\xymatrix{ {\phi} \ar@{-}[dd] \ar[d] \ar@{-}[rr]^(0.6){H} \ar[r]
& & {\psi}\ar@{-}[dd] \ar[d]\\
&&\\
 {g \phi g^{-1}} \ar@{-}[rr]^(0.6){g H g^{-1}} \ar[r]& & {g \phi g^{-1}} }
$$
As the components are contractible, we can fill in the rectangle
whose boundary is made out of homotopies. This proves that $H$
is $G$-equivariant up to homotopy.
\end{proof}

\begin{prop}\label{prop:Harrison}
The Harrison set $\Har(A,G)$  is natural in $G$. When $G$ is
an abelian group, then $\Har(A,G)$ is an abelian group as well.
\end{prop}
\begin{proof}
We give a translation of Greither's proof from~\cite{Greither:LNM1534}
to the topological setting. Given a group homomorphism $\phi\:G\lra H$,
for a $G$-Galois extension $B/A$ we define
\[
\phi_*B=F(H_+,B)^{\mathrm{h}G} = F(EG_+,F(H_+,B))^G,
\]
the homotopy fixed points of the function spectrum $F(H_+,B)$ with
the action coming from the left $G$-action on $H$ and $B$. Then
$\phi_*B$ has a natural $H$-action induced from the right action
of~$H$ on itself.

Choose a  free contractible left $G$-space $EG$ and a free contractible
left $H$-space $EH$ and also write $E'H$ for $EH$ with the trivial
$H$-action but also viewed as a left $G$-space through the homomorphism
$\phi$. We also view $H$ as a left $G$-space via the homomorphism $\phi$
and as a left $H$-space by inverse right multiplication. There is a left
$G$-action on $H\times_H EH$ induced from the action on the $H$ factor
and the isomorphism
\[
H\times_H EH\iso E'H; \quad [\eta,x]\longleftrightarrow \eta x
\]
is $G$-equivariant.

Consider a chain of homomorphisms of groups $G
\xrightarrow{\phi}H\xrightarrow{\psi}K$. We have to prove that
\begin{equation}\label{eq:natural}
\psi_*(\phi_*B) \sim (\psi\o\phi)_*B.
\end{equation}
To this end we have to compare
$F(K_+,F(H_+,B)^{\mathrm{h}G})^{\mathrm{h}H}$ with
$F(K'_+,B)^{\mathrm{h}G}$, where $K'$ denotes $K$ with the left
action of $G$ via $\psi\o\phi$. There is a chain of equivalences
\begin{align*}
F(K_+,F(H_+,B)^{\mathrm{h}G})^{\mathrm{h}H} &=
 F(EH_+,F(K_+,F(EG_+,F(H_+,B))^G))^H \\
&\iso
 F(EG\times (H \times EH)\times_H K_+,B)^G \\
&\iso
 F(EG\times E'H \times K'_+,B)^G \\
&\stackrel{\sim}{\lla} F(EG \times K'_+,B)^G \\
&\iso F(EG_+,F(K'_+,B))^G=F(K'_+,B)^{\mathrm{h}G}.
\end{align*}

We consider $B \wedge_A \phi_*B$. The left-hand factor of $B$
has a trivial $G$-action and  $B$ is strongly dualizable. Hence
it is equivalent to $\phi_*(B \wedge_A B)$ by the following chain
of identifications:
\begin{align*}
B \wedge_A \phi_*B & = B \wedge_A F(EG_+,F(H_+,B))^G\\
                   & \cong B \wedge_A F((EG \times H)_+,B)^G\\
                  & \stackrel{\sim}{\lra} F((EG \times H)_+, B \wedge_A B)^G\\
                   & \cong  F(EG_+, F(H_+, B \wedge_A B))^G\\
                   & = \phi_*(B \wedge_A B).
\end{align*}
Using Proposition~\ref{prop:FFbaseChange} (c) (or results
of~\cite{R:Opusmagnus}) it is enough to check that $B \wedge_A \phi_*B$
is an $H$-Galois extension of $B$.

In order to check this, we consider the two natural inclusions
of the trivial group
\[
G \stackrel{i}{\hookleftarrow} e \stackrel{j}{\hookrightarrow} H.
\]
It is obvious that $i_*B \cong \prod_G B$ and $j_*B \cong \prod_H B$.
Using naturality~\eqref{eq:natural} and $\phi \circ i =j$, we obtain
\[
\phi_*(B\wedge_A B) \cong \phi_*(\prod_G B) \cong \phi_*(i_* B)
\stackrel{\sim}{\lla} j_*B \cong \prod_H B.
\]
Therefore $\phi_*B$ is unramified with respect to the $H$-action
and its $H$-homotopy fixed points agree with $A$.

If we consider abelian Galois extensions, then the source and
target in~\eqref{eqn:Harrison-Natural} have abelian group structures.
On the algebraic side, the map induced by the abelian multiplication
$\mu\:G\times G\lra G$ is a homomorphism which sends two $G$-Galois
extensions $B_*/A_*$ and $C_*/A_*$ to
$\mu_*(B_*\otimes_{A_*}C_*)/A_*$. Since the multiplication
homomorphism is surjective, there is a short exact sequence
\[
0\ra K=\ker\mu\lra G \times G \xrightarrow{\mu}G\ra 0
\]
and so using Harrison's formula~\cite[p.3]{Harrison:AbExtns}
we obtain
\[
\mu_*(B_* \otimes_{A_*} C_*) = (B_* \otimes_{A_*} C_*)^K.
\]
Mimicking this in the geometric situation we set
$B\.C=(B\wedge_A C)^{\mathrm{h}K}$ for any abelian $G$-Galois
extensions $B$ and $C$ of $A$. The proof that the induced map
$\phi_*\:\Har(A,G) \lra \Har(A,H)$ of every homomorphism between
abelian groups $\phi\: G \lra H$ is a homomorphism only uses
naturality and is analogous to the
proof of~\cite[3.2]{Greither:LNM1534}.
\end{proof}
\begin{rem}\label{rem:Harrison-Naturality}
Let $A$ be a commutative $\mathbb{S}$-algebra. \\
(a) The Harrison functor $\Har(A,-)$ restricted to abelian
groups is additive: for abelian groups $G_1$ and $G_2$,
\begin{align*}
\Har(A,G_1\times G_2)&\iso\Har(A,G_1)\times\Har(A,G_2).
\end{align*} \\
(b) If $G$ is an abelian group of exponent $n$, then multiplication
by~$n$ is induced by multiplication by~$n$ on $G$ which factors
through the trivial group, so
\[
n\Har(A,G) = 0.
\]
\end{rem}

It would be interesting to have a better understanding of the
function out of $\Har(A,G)$ into a subcategory of the category
of $A_*$-algebras which sends a $G$-Galois extension $B/A$ to
$B_*/A_*$. In Section~\ref{sec:TopAbelExtns&KummerThy} we
investigate the corresponding question for $\Har^w(A,G)$ in
the case where $G$ is abelian.

\section{Topological Kummer theory}\label{sec:TopAbelExtns&KummerThy}

We will now describe analogous constructions to those of
Section~\ref{sec:Abel&Kummer} when the following conditions
are satisfied.
\begin{cond}\label{cond:A-G}
$A$ is a commutative $\mathbb{S}$-algebra and $G$ is a finite
abelian group for which
\begin{itemize}
\item
$A_0$ contains $1/|G|$;
\item
$A_0$ contains a primitive $d$-th root of unity $\zeta$,
where $d$ is the exponent of $G$;
\item
$A_0$ is connected (\ie, it has no non-trivial idempotents).
\end{itemize}
\end{cond}

By generalizing constructions of~\cite{SVW} as in
Example~\ref{ex:Cyclotomic} we can always arrange for the second
condition to hold whenever the first does. The third condition
is not strictly necessary but simplifies the ensuing discussion.
\begin{thm}\label{thm:Kumm-Top->Alg}
Suppose that $A$ and $G$ satisfy \emph{Condition~\ref{cond:A-G}}
and let $B/A$ be a $G$-Galois extension. If for every invertible
$A$-module $U$, $U_*$ is an invertible graded $A_*$-module, then
$B_*/A_*$ is a $G$-Galois extension.
\end{thm}

We will see later, that the invertibility condition is not void
in general. For a more thorough treatment of the question of when
invertible module spectra give rise to invertible coefficients,
see~\cite{AB&BR:Inv}.
\begin{proof}
Notice that there is a decomposition of the form~\eqref{eqn:R[G]-Decomp},
\[
A_0[G]=\bigoplus_\chi A_0(\chi)
\]
defined using idempotents $e_\chi$ as defined
in~\eqref{eqn:R[G]-Idempot}. By
Theorem~\ref{thm:BraveNewGalois-Properties}(c),
$B\wedge G_+\sim F_A(B,B)$. Each $e_\chi$ is an element of $A_0[G]$
and can be realized by a map of $A$-modules $A \ra A \wedge G_+$.
Composing this with the unit of $B$ gives rise to a map
\[
A \lra B\wedge G_+\sim F_A(B,B),
\]
whose action on $B$ can be iterated to produce an $A$-module
$B(\chi)=e_\chi B$ which is well defined up to homotopy equivalence.
There is a homotopy decomposition of $A$-modules
\begin{equation}\label{eqn:Decomp-B}
B\sim\bigvee_{\chi}B(\chi).
\end{equation}
As in the algebraic case, there are also pairings
\begin{equation}\label{eqn:B(chi)-Product}
B(\chi_1)\Smash_AB(\chi_2)\lra B(\chi_1\chi_2).
\end{equation}
Now smash with $B$ and recall that $B$ is a faithful $A$-module.
We can consider the extension $B\Smash_AB/B$ which is equivalent
to $(\prod_{\gamma\in G}B)/B$. The decomposition analogous to
that of~\eqref{eqn:Decomp-B},
\begin{equation*}
B\Smash_AB\sim\bigvee_{\chi}(B\Smash_AB)(\chi),
\end{equation*}
is induced from that of $B$ by smashing with $B$ and
\[
(B\Smash_AB)(\chi)\sim B\Smash_A(B(\chi)).
\]
For the product maps we also have homotopy commutative diagrams
\[
\xymatrix{ {B\Smash_A(B(\chi_1)\Smash_A B(\chi_2))} \ar[r] \ar[d]
                                & {B\Smash_AB(\chi_1\chi_2)} \ar[d] \\
{(B\Smash_AB(\chi_1))\Smash_B(B\Smash_AB(\chi_2))}\ar[r] &
{(B\Smash_AB)(\chi_1\chi_2)} }
\]
so if we can show that the bottom maps are equivalences of
$B$-modules then since $B$ is a faithful $A$-module, the maps
of~\eqref{eqn:B(chi)-Product} are equivalences of $A$-modules.
But the necessary verification is formally similar to that for
the algebraic case proved in~\cite{Greither:LNM1534} since in
homotopy there is an isomorphism of $A_*[G]$-modules
\[
\pi_*(B\Smash_AB)\iso\prod_{\gamma\in G}B_*.
\]

Thus we have shown that the map of~\eqref{eqn:B(chi)-Product}
is a weak equivalence. In particular, each $B(\chi)$ is an
invertible $A$-module. Now by assumption $B(\chi)_*$ is an
invertible graded $A_*$-module and so is projective. From
this we conclude that $B_*$ is a direct sum of projective
$A_*$-modules and so the K\"unneth spectral sequence collapses
to give
\[
B^A_*B\iso B_*\oTimes_{A_*}B_*,
\]
and therefore
\[
B_*\oTimes_{A_*}B_*=\prod_{\gamma\in G}B_*.
\]
As the order of the group is inverted in $A_0$ we also have
\[
(B_*)^G={(B^{\mathrm{h}G})}_*=A_*,
\]
and therefore $B_*/A_*$ is a $G$-Galois extension.
\end{proof}

Here is a reinterpretation of what we have established by combining
Theorems~\ref{thm:Obstructions} and~\ref{thm:Kumm-Top->Alg}.
\begin{thm}\label{thm:Kumm-Top=Alg}
Suppose that $A$ and $G$ satisfy \emph{Condition~\ref{cond:A-G}}
and that the coefficients of invertible $A$-modules are invertible
graded $A_*$-modules. Then there is a natural bijection of sets
\[
\Real_{A,G}\:\Har(A_*,G)\xrightarrow{\iso}\Har^w(A,G).
\]
Therefore the weak Harrison set is actually a group, because it
inherits the group structure from the algebraic Harrison group.
\end{thm}
\begin{rem}
Note that under the above assumptions taking homotopy groups
always gives a map from either Harrison set to $\Har(A_*,G)$,
but for $\Har(A,G)$ we do not obtain an isomorphism, because
we do not have an inverse  map.
\end{rem}

We close this section with an example where one can classify all
of the topological Kummer extensions of an $\mathbb{S}$-algebra.
\begin{ex}\label{ex:Har-KO}
We have
\[
\Har^w(KO[1/2],C_2)\iso C_2\times C_2\times C_2.
\]
To see this, first from~\eqref{eqn:KO*} we see that
\[
KO[1/2]_*=\Z[1/2][y,y^{-1}],
\]
where $y\in KO[1/2]_4$. As we established in~\cite{AB&BR:Inv} that
invertible $KO[1/2]$-modules have invertible coefficients,
Theorem~\ref{thm:Kumm-Top=Alg} yields
\[
\Har^w(KO[1/2],C_2)\iso\Har(\Z[1/2][y,y^{-1}],C_2).
\]
By Proposition~\ref{prop:KummSeq}, we find that
\begin{align*}
\Har^w(KO[1/2],C_2)
&\iso(\Z[1/2][y,y^{-1}])^\times/\left((\Z[1/2][y,y^{-1}])^\times\right)^2 \\
&\iso C_2\times C_2\times C_2,
\end{align*}
with generators the cosets of $-1,2,y$ with respect to
$\left((\Z[1/2][y,y^{-1}])^\times\right)^2$.

This leads to three non-trivial $C_2$-extensions of $KO[1/2]$ with
coefficient rings
\[
KO_*[1/2,i],\; KO_*[1/2,\sqrt{2}],\; KO_*[1/2,i\sqrt{2}]
\]
which correspond to the cosets of $-1,2,-2$ and the arithmetic
extensions
\[
\Z[1/2,i]/\Z[1/2],\; \Z[1/2,\sqrt{2}]/\Z[1/2],\;
\Z[1/2,i\sqrt{2}]/\Z[1/2].
\]
The extension $KU[1/2]/KO[1/2]$ corresponds to
$KO_*[1/2,\sqrt{y/2}]$. There are three  more exotic
extensions associated with the rings $KO_*[1/2,\sqrt{y}]$,
$KO_*[1/2,i\sqrt{2y}]$ and $KO_*[1/2,i\sqrt{y}]$. These
are the Harrison products of  $KU[1/2]$ with the three
above.

Note that we can adjoin an $8$-th root of unity to $KO[1/2]$
and get a $C_2 \times C_2$-extension (see Example~\ref{ex:Cyclotomic}).
By taking homotopy fixed points with respect to subgroups
we obtain the arithmetic $C_2$-extensions listed above.
$$
\xymatrix{
{}& {KO[1/2,\zeta_8]} \ar@{-}[dl]\ar@{-}[d]\ar@{-}[dr] & {} \\
{KO[1/2,i]} \ar@{-}[dr] & {KO[1/2,\sqrt{2}]} \ar@{-}[d] &
{KO[1/2,i\sqrt{2}]} \ar@{-}[dl] \\
{}& {KO[1/2]} & {} }
$$
\end{ex}

\section{Topological abelian extensions}\label{sec:TopAbelianExtns}

In the following we will consider finite abelian extensions $B/A$
without assuming that the order of the Galois group $G$ is
invertible in $A_*$.
\begin{thm}\label{thm:GalExtn-InvtbleA[G]}
For every finite abelian $G$-Galois extension $B/A$, $B$ is an
invertible $A[G]$-module.
\end{thm}
\begin{proof}
For $G$ abelian we have a natural evaluation map
\[
\epsilon\: F_{A[G]}(B, A[G]) \wedge_{A[G]} B \lra A[G].
\]
We will prove that this map is an equivalence. As $B$ is a faithful
$A$-module, it suffices to consider the map $B \wedge_A \epsilon$
instead. As $B$ is dualizable over $A$ and self-dual, we can
identify $B \wedge_A F_{A[G]}(B, A[G]) \wedge_{A[G]} B$ with
$F_{A[G]}(B, B[G]) \wedge_{A[G]} B$. Inducing up to $B$ then yields
an equivalence with $F_{B[G]}(B \wedge_A B, B[G]) \wedge_{B[G]} B
\wedge_A B$. Then $B \wedge_A \epsilon$ factors as in the following
diagram.
\[
\xymatrix{ {B \wedge_A F_{A[G]}(B, A[G]) \wedge_{A[G]} B}
\ar[d]_{\sim}
\ar[rrr]^{B \wedge_A \epsilon} &&& {B \wedge_A A[G]} \ar[d]^{\iso}\\
{F_{A[G]}(B, B[G]) \wedge_{A[G]} B} \ar[rrr] \ar[d]_{\iso}&
& &  {B[G]} \\
{F_{B[G]}(B \wedge_A B, B[G]) \wedge_{B[G]} B \wedge_A B}
\ar[d]_{\sim}^{\Theta} & & &\\
{F_{B[G]}(B \wedge_A B, B[G]) \wedge_{B[G]} F(G_+,B)}
& & &\\
{F_{B[G]}(B[G], B[G]) \wedge_{B[G]} B[G]} \ar[u]^{\sim}_{\Upsilon}
\ar@/_6ex/[uuurrr]_(0.6){\iso} & & & }
\]
Here we use the equivalence
$$\begin{CD} B \wedge_A B @>\sim>\Theta>
F(G_+,B) @<\sim<\Upsilon< B[G]
\end{CD}$$
where $\Upsilon$ is the topological analogue of~\eqref{eqn:groupring},
in particular it is an equivalence of $B[G]$-modules.
\end{proof}
\begin{ex}\label{ex:Pic-KU}
The result of Theorem~\ref{thm:GalExtn-InvtbleA[G]} gives rise
to examples of invertible $A[G]$-modules whose coefficient groups
do not yield elements in the algebraic Picard group $\Pic(A_*[G])$.
Consider for instance the $C_2$-Galois extension $KU/KO$. We
know that $KU$ is an invertible $KO[C_2]$-module, but $KU_*$
is definitely not an invertible $KO_*[C_2]$-module, because
it is not even projective.
\end{ex}

The Harrison group $\Har(A,G)$ is related to the Picard group
of the group ring $A[G]$. We will make use of the constructions
introduced in Proposition~\ref{prop:Harrison}.
\begin{thm}\label{thm:Har->Pic}
There is a homomorphism
\[
\Psi_G\:\Har(A,G) \lra \Pic(A[G]).
\]
In particular,  for every finite abelian group\/ $G$ of exponent~$n$,
the image of $\Psi_G$ is contained in the $n$-torsion subgroup
of $\Pic(A[G])$.
\end{thm}
\begin{proof}
For a finite abelian group $G$ we define
\[
\Psi_G\:\Har(A,G) \lra \Pic(A[G]); \quad \Psi_G([B]) = [B],
\]
where the first equivalence class is in the Harrison group of
$G$-Galois extensions of $A$ and the second denotes an isomorphism
class in the homotopy category of $A[G]$-modules. Whenever we have
to choose a representing element $B'$ for $[B]$ in $\Pic(A[G])$ it
will be a cofibrant $A[G]$-module. Equivalent $G$-Galois extensions
over $A$ are in particular equivalent $A$-modules with $G$-action;
therefore $\Psi_G$ is well-defined.

We have to show that it is a homomorphism, \ie, that the isomorphism
class of $\mu_*(B' \wedge_A C') = (B' \wedge_A C')^{\mathrm{h}K}$
coincides with that of $B' \wedge_{A[G]} C'$. Here
$\mu\: G \times G \lra G$ is the multiplication in $G$ and $K=\ker\mu$,
while $B'$ and $C'$ are $A[G]$-cofibrant models of $B$ and $C$
respectively.

If $M$ is a cofibrant $A$-module, then $EG_+\wedge M$ is a cofibrant
model for $M$ in the category of $A[G]$-modules: the fibrations
and weak equivalences in the categories $\mathscr{M}_A$ and
$\mathscr{M}_{A[G]}$ are defined via the forgetful functor to
$\mathscr{M}_\mathbb{S}$. Therefore if $U$ denotes the forgetful
functor from $\mathscr{M}_{A[G]}$ to $\mathscr{M}_A$, the lifting
diagram
\[
\xymatrix{
{*} \ar[r] \ar[d]
& {U(X)} \ar[d]^{U(f)} \\
{M} \ar@{.>}[ur]^{\xi} \ar[r]  & {U(Y)}
}
\]
for any acyclic fibration $f\: X \lra Y$ in $\mathscr{M}_{A[G]}$
has a $G$-equivariant extension
\[
\xymatrix{
{*} \ar[r] \ar[d]
& {X} \ar[d]^{f} \\
{EG_+ \wedge M} \ar@{.>}[ur]^{\bar{\xi}} \ar[r] & {Y}
}
\]
and therefore we can identify $B' \wedge_{A[G]} C'$ with
$EG_+ \wedge B\wedge_{A[G]} EG_+ \wedge C$ and this in turn
is equivalent to $(E(G \times G)_+ \wedge B \wedge_A C)/K$.
We are left with the identification of the homotopy orbits
$(E(G\times G)_+\wedge B\wedge_A C)/K=(B\wedge_A C)_{\mathrm{h}K}$
and the homotopy fixed points.

As $B$ and $C$ are both dualizable over $A$, we obtain the chain
of identifications
\begin{align*}
(B \wedge_A C) \wedge_A (B \wedge_A C)_{\mathrm{h}K} \sim
&(B \wedge_A C\wedge_A B \wedge_A C)_{\mathrm{h}K}
\sim (\prod_{G \times G} (B \wedge_A C))_{\mathrm{h}K} \\
\sim & (\prod_{G \times G} (B \wedge_AC))^{\mathrm{h}K} \sim
(B \wedge_A C) \wedge_A (B \wedge_A C)^{\mathrm{h}K},
\end{align*}
and this shows that
$(B\wedge_A C)_{\mathrm{h}K} \sim (B\wedge_A C)^{\mathrm{h}K}$,
since $B \wedge_A C$ is faithful over $A$. Here we use the fact
that $\prod_{G \times G}(B \wedge_A C)$ is equivalent to a wedge
of copies of $F(K_+, B \wedge_A C)$ which is a free $K$-spectrum,
and hence it has a trivial Tate spectrum~\cite[proposition~2.4]{JG&JM:Tate}.

By Remark~\ref{rem:Harrison-Naturality}(b), when $G$ is an abelian
group of exponent~$n$, the image of $\Psi_G$ is contained in the
$n$-torsion subgroup of $\Pic(A[G])$.
\end{proof}

\section{Units of Galois extensions}\label{sec:Units}

One instance of Hilbert's theorem~90 involves the vanishing of the
first cohomology of Galois groups with coefficients in the units of
a field extension and the straightforward identification of the
fixed points of the units of the extension with the units in the
base field. Note that for a general $G$-Galois extension of rings
$S/R$ Hilbert's theorem~90 does not hold: instead of a vanishing
result for the first group cohomology, there is an exact
sequence~\cite[\S5]{CHR},
\begin{multline*}
\ph{\Br}
0\ra\mathrm{H}^1(G,S^\times)\lra\Pic(R)\lra\mathrm{H}^0(G,\Pic(S))
\lra\mathrm{H}^2(G,S^\times)        \\
\lra\Br(S/R)\lra\mathrm{H}^1(G,\Pic(S))\lra\mathrm{H}^3(G,S^\times),
\ph{\Br}
\end{multline*}
in which $\Br(S/R)$ is the relative Brauer group. Examples of Galois
extensions with non-trivial $\mathrm{H}^1(G,S^\times)$ are mentioned
in~\cite[5.5~(d)]{CHR}. We will prove a version of Hilbert's
theorem~90 corresponding to the classical statement for invariants.

In the following we use the concept of units of ring spectra. These
were introduced by Patterson, Stong and Waldhausen and their
multiplicative properties were developed in~\cite{JPM:Einftyrings}.
More material on these is contained in the notes~\cite{Ando} and we
are grateful to M.~Ando for providing them.
\begin{defn}\label{defn:GL1}
Let $R$ be a ring spectrum. The units $GL_1(R)$ of the ring
spectrum~$R$ are defined via the following homotopy pullback square:
\[
\xymatrix{
{GL_1(R)} \ar[r] \ar[d] & {\Omega^\infty(R)} \ar[d] \\
{(\pi_0(R))^\times} \ar[r] & {\pi_0(R)}}
\]
in which $\Omega^\infty(R)$ denotes the underlying infinite loop
space of the $\Omega$-spectrum associated to~$R$, \ie, the zeroth
space of the spectrum~$R$.
\end{defn}
The quotient map from the space $\Omega^\infty(R)$ to its path
components  is a fibration. Therefore the units of $R$ are given by
an actual pullback square.

Now assume that $R$ possesses an action of some finite group~$G$ by
maps of ring spectra. More precisely, let $R$ be a naive
$G$-spectrum with a coherent $E_\infty$-structure in the sense
of~\cite[VII Def.~2.1]{LMS}. We recall from~\cite[VII
proposition~2.8]{LMS} that the zeroth space $R(0)=\Omega^\infty(R)$
inherits a $G$-$E_\infty$-structure from~$R$.

The homotopy groups of $R$ inherit the $G$-action as well. As
everything takes place in a setting of $\Omega$-spectra, the zeroth
homotopy group $\pi_0(R)$ is given by $\pi_0(\Omega^\infty(R))$. The
inclusion of the units into the full ring $\pi_0(R)$ is clearly
$G$-equivariant and so is the quotient map from $\Omega^\infty(R)$
to $\pi_0(R)$.
\begin{thm}\label{thm:GL1(R)hG}
Assume that $R$ is a $G$-ring spectrum as above and for which
$\pi_0(R^{\mathrm{h}G})\cong(\pi_0(R))^G$. Then the homotopy fixed
points of the units $GL_1(R)$ are given by the units of
$R^{\mathrm{h}G}$.
\end{thm}
\begin{proof}
Taking the zeroth space of a spectrum commutes with homotopy fixed
points, because using the setting of~\cite[I,\S 3]{LMS} we have the
following chain of identifications
\begin{align*}
(\Omega^\infty(R))^{\mathrm{h}G} &= F( EG_+,\Omega^\infty(R))^G
= F( EG_+, R(0))^G \\
&= (F( EG_+, R)(0))^G
= (F( EG_+, R)^G)(0) \\
&=\Omega^\infty(R^{\mathrm{h}G}).
\end{align*}
By assumption, the homotopy fixed points of the discrete set
$\pi_0(R)$ are
\[
F( EG_+,\pi_0(R))^G\cong (\pi_0(R))^G \cong\pi_0(F( EG_+, R)^G)
\]
and therefore the pullback for the homotopy fixed points of
$GL_1(R)$ is the pullback of the diagram
\[
\xymatrix{
{}& {\Omega^\infty(R^{\mathrm{h}G})} \ar[d] \\
{(\pi_0(R^{\mathrm{h}G}))^\times}\ar[r] & {\pi_0(R^{\mathrm{h}G})} }
\]
and by definition this is $GL_1(R^{\mathrm{h}G})$.
\end{proof}

We now obtain a topological version of Hilbert's theorem~90 as an
immediate consequence of the above result. Note that the following
result also holds if $B$ is ramified over~$A$.
\begin{cor}\label{cor:Hilbert90}
Let $G$ be a finite group and let $A\lra B$ be a weak $G$-Galois
extension with $(\pi_0(B))^G\cong\pi_0(A)$. Then
\[
GL_1(B)^{\mathrm{h}G}\sim GL_1(A).
\]
\end{cor}

The condition on the zeroth homotopy group is of course satisfied
in the case of a realization of an  algebraic  $G$-Galois extension
of a commutative $\mathbb{S}$-algebra. In the special case of
Eilenberg-Mac~Lane spectra $HR\lra HS$ the result gives the
classical identity $(S^\times)^G\cong R^\times$. The $C_2$-Galois
extension $KU/KO$ also satisfies the condition on $\pi_0$ and so
do the naturally occurring examples in~\cite[\S 5]{R:Opusmagnus}.
However, there are examples where this condition is not satisfied:
take $A$ to be $\bigvee_{n \in \Z} \Sigma^{2n} KO$ and take $B$
to be
$$
\bigvee_{n \in \Z} \Sigma^{2n} KO \wedge_{KO} KU \sim \bigvee_{n \in \Z}
\Sigma^{2n} KU.
$$
Therefore $\pi_0(B)^{C_2}$ consists of copies of the integers, whereas
$\pi_0(A)$ contains summands $\Z/2\Z$.

\end{document}